\documentclass[final,3p]{elsarticle}
 \usepackage{graphics}
 \usepackage{graphicx}
 \usepackage{epsfig}
\usepackage{amssymb}
 \usepackage{amsthm}
 \usepackage{lineno}
 \usepackage{amsmath}
   \numberwithin{equation}{section}
\usepackage{mathrsfs}

\newtheorem{thm}{Theorem}[section]

\newtheorem{lem}[thm]{Lemma}

\biboptions{sort&compress,square}
\journal{some mathematical journals}
\begin{document}
\begin{frontmatter}
\author[rvt1]{Sining Wei}
\ead{weisn835@nenu.edu.cn}
\author[rvt2]{Jian Wang}
\ead{wangj484@nenu.edu.cn}
\author[rvt3]{Yong Wang\corref{cor2}}
\ead{wangy581@nenu.edu.cn}
\cortext[cor2]{Corresponding author.}

\address[rvt1]{School of Data Science and Artificial Intelligence, Dongbei University of Finance and Economics, \\
Dalian, 116025, P.R.China}
\address[rvt2]{School of Science, Tianjin University of Technology and Education, Tianjin, 300222, P.R.China}
\address[rvt3]{School of Mathematics and Statistics, Northeast Normal University, Changchun, 130024, P.R.China}

\title{ Conformal Perturbations
of Twisted Dirac Operators and \\Noncommutative residue }
\begin{abstract}
In this paper, we obtain two kinds of Kastler-Kalau-Walze type theorems
for conformal perturbations
of twisted Dirac operators and conformal perturbations
of signature operators by a vector bundle with a non-unitary connection on six-dimensional manifolds with  (respectively without)boundary.
\end{abstract}
\end{frontmatter}
\section{Introduction}
\label{1}
The noncommutative residue found in \cite{Gu,Wo} plays a prominent role in noncommutative geometry. For one-dimensional manifolds,
the noncommutative residue was discovered by Adler \cite{MA} in connection with geometric aspects of nonlinear partial differential equations. For arbitrary closed compact $n$-dimensional manifolds, the noncommutative residue was introduced by Wodzicki in \cite{Wo} using the theory of zeta functions of elliptic pseudodifferential operators. In \cite{Co1}, Connes used the noncommutative residue to derive a conformal 4-dimensional Polyakov action analogy.
Furthermore, Connes made a challenging observation that the noncommutative residue of the square of the inverse of the
Dirac operator was proportional to the Einstein-Hilbert action in \cite{Co2}.
In \cite{Ka}, Kastler gave a brute-force proof of this theorem. In \cite{KW}, Kalau and Walze proved this theorem in the
normal coordinates system simultaneously. And then, Ackermann proved that
the Wodzicki residue  of the square of the inverse of the
Dirac operator ${\rm  Wres}(D^{-2})$ in turn is essentially the second coefficient
of the heat kernel expansion of $D^{2}$ in \cite{Ac}.

In \cite{RP}, Ponge defined lower dimensional volumes of Riemannian manifolds by the Wodzicki residue. Fedosov et al. defined a noncommutative residue on Boutet de Monvel's algebra and proved that it was a unique continuous trace in \cite{FGLS}. In \cite{S}, Schrohe gave the relation between the Dixmier trace and the noncommutative residue for manifolds with boundary. In \cite{Wa3}, Wang generalized the Kastler-Kalau-Walze type theorem to the cases of 3, 4-dimensional spin manifolds with boundary and proved a Kastler-Kalau-Walze type theorem.  In \cite{Wa3,Wa4,WJ2,WJ3,WJ4}, Y. Wang and his coauthors computed the lower dimensional volumes for 5, 6, 7-dimensional spin manifolds with boundary and also got some Kastler-Kalau-Walze type theorems. In \cite{WJ5}, authors computed $\widetilde{{\rm Wres}}[(\pi^+D^{-2})\circ(\pi^+D^{-n+2})]$ for any-dimensional manifolds with boundary, and proved a general Kastler-Kalau-Walze type theorem.

In \cite{WJ6}, J. Wang and Y. Wang proved two kinds of Kastler-Kalau-Walze type theorems
for conformal perturbations
of twisted Dirac operators and conformal perturbations
of signature operators by a vector bundle with a non-unitary connection on four-dimensional manifolds with  (respectively without)boundary.

{\bf The motivation of this paper} is to establish two Kastler-Kalau-Walze type theorems for conformal perturbations
of twisted Dirac operators and conformal perturbations
of signature operators with non-unitary connections on six-dimensional manifolds with boundary. We know that the leading symbol of conformal perturbations
of twisted Dirac operators is not $\sqrt{-1}c(\xi)$. This is the reason that we study the residue of conformal perturbations
of twisted Dirac operators.\\

 This paper is organized as follows: In Section 2, we recall
some basic facts and formulas about Boutet de Monvel's calculus. In Section 3, we give a Kastler-Kalau-Walze type theorems for conformal perturbations
of twisted Dirac operators on six-dimensional manifolds with boundary. In Section 4 and Section 5, we recall the definition of conformal perturbations
of signature operators and compute their symbols, and we give a Kastler-Kalau-Walze type theorems for conformal perturbations
of signature operators on six-dimensional manifolds with boundary.

\section{Boutet de Monvel's calculus and noncommutative residue}
In this section, we shall recall
some basic facts and formulas about Boutet de Monvel's calculus.
Let $$ F:L^2({\bf R}_t)\rightarrow L^2({\bf R}_v);~F(u)(v)=\int e^{-ivt}u(t)\texttt{d}t$$ denote the Fourier transformation and
$\varphi(\overline{{\bf R}^+}) =r^+\varphi({\bf R})$ (similarly define $\varphi(\overline{{\bf R}^-}$)), where $\varphi({\bf R})$
denotes the Schwartz space and
  \begin{equation}
r^{+}:C^\infty ({\bf R})\rightarrow C^\infty (\overline{{\bf R}^+});~ f\rightarrow f|\overline{{\bf R}^+};~
 \overline{{\bf R}^+}=\{x\geq0;x\in {\bf R}\}.
\end{equation}
We define $H^+=F(\varphi(\overline{{\bf R}^+}));~ H^-_0=F(\varphi(\overline{{\bf R}^-}))$ which are orthogonal to each other. We have the following
 property: $h\in H^+~(H^-_0)$ if and only if $h\in C^\infty({\bf R})$ which has an analytic extension to the lower (upper) complex
half-plane $\{{\rm Im}\xi<0\}~(\{{\rm Im}\xi>0\})$ such that for all nonnegative integer $l$,
 \begin{equation}
\frac{d^{l}h}{d\xi^l}(\xi)\sim\sum^{\infty}_{k=1}\frac{d^l}{d\xi^l}(\frac{c_k}{\xi^k})
\end{equation}
as $|\xi|\rightarrow +\infty,{\rm Im}\xi\leq0~({\rm Im}\xi\geq0)$.

 Let $H'$ be the space of all polynomials and $H^-=H^-_0\bigoplus H';~H=H^+\bigoplus H^-.$ Denote by $\pi^+~(\pi^-)$ respectively the
 projection on $H^+~(H^-)$. For calculations, we take $H=\widetilde H=\{$rational functions having no poles on the real axis$\}$ ($\tilde{H}$
 is a dense set in the topology of $H$). Then on $\tilde{H}$,
 \begin{equation}
\pi^+h(\xi_0)=\frac{1}{2\pi i}\lim_{u\rightarrow 0^{-}}\int_{\Gamma^+}\frac{h(\xi)}{\xi_0+iu-\xi}d\xi,
\end{equation}
where $\Gamma^+$ is a Jordan close curve included ${\rm Im}\xi>0$ surrounding all the singularities of $h$ in the upper half-plane and
$\xi_0\in {\bf R}$. Similarly, define $\pi^{'}$ on $\tilde{H}$,
 \begin{equation}
\pi'h=\frac{1}{2\pi}\int_{\Gamma^+}h(\xi)d\xi.
\end{equation}
So, $\pi'(H^-)=0$. For $h\in H\bigcap L^1(R)$, $\pi'h=\frac{1}{2\pi}\int_{R}h(v)dv$ and for $h\in H^+\bigcap L^1(R)$, $\pi'h=0$.
Denote by $\mathcal{B}$ Boutet de Monvel's algebra (for more details, see Section 2 of \cite{Wa5}).

An operator of order $m\in {\bf Z}$ and type $d$ is a matrix
$$A=\left(\begin{array}{lcr}
  \pi^+P+G  & K  \\
   T  &  S
\end{array}\right):
\begin{array}{cc}
\   C^{\infty}(X,E_1)\\
 \   \bigoplus\\
 \   C^{\infty}(\partial{X},F_1)
\end{array}
\longrightarrow
\begin{array}{cc}
\   C^{\infty}(X,E_2)\\
\   \bigoplus\\
 \   C^{\infty}(\partial{X},F_2)
\end{array},
$$
where $X$ is a manifold with boundary $\partial X$ and
$E_1,E_2~(F_1,F_2)$ are vector bundles over $X~(\partial X
)$.~Here,~$P:C^{\infty}_0(\Omega,\overline {E_1})\rightarrow
C^{\infty}(\Omega,\overline {E_2})$ is a classical
pseudodifferential operator of order $m$ on $\Omega$, where
$\Omega$ is an open neighborhood of $X$ and
$\overline{E_i}|X=E_i~(i=1,2)$. $P$ has an extension:
$~{\cal{E'}}(\Omega,\overline {E_1})\rightarrow
{\cal{D'}}(\Omega,\overline {E_2})$, where
${\cal{E'}}(\Omega,\overline {E_1})~({\cal{D'}}(\Omega,\overline
{E_2}))$ is the dual space of $C^{\infty}(\Omega,\overline
{E_1})~(C^{\infty}_0(\Omega,\overline {E_2}))$. Let
$e^+:C^{\infty}(X,{E_1})\rightarrow{\cal{E'}}(\Omega,\overline
{E_1})$ denote extension by zero from $X$ to $\Omega$ and
$r^+:{\cal{D'}}(\Omega,\overline{E_2})\rightarrow
{\cal{D'}}(\Omega, {E_2})$ denote the restriction from $\Omega$ to
$X$, then define
$$\pi^+P=r^+Pe^+:C^{\infty}(X,{E_1})\rightarrow {\cal{D'}}(\Omega,
{E_2}).$$
In addition, $P$ is supposed to have the
transmission property; this means that, for all $j,k,\alpha$, the
homogeneous component $p_j$ of order $j$ in the asymptotic
expansion of the
symbol $p$ of $P$ in local coordinates near the boundary satisfies:
$$\partial^k_{x_n}\partial^\alpha_{\xi'}p_j(x',0,0,+1)=
(-1)^{j-|\alpha|}\partial^k_{x_n}\partial^\alpha_{\xi'}p_j(x',0,0,-1),$$
then $\pi^+P:C^{\infty}(X,{E_1})\rightarrow C^{\infty}(X,{E_2})$
by Section 2.1 of \cite{Wa5}.

In the following, write $\pi^+D^{-1}=\left(\begin{array}{lcr}
  \pi^+D^{-1}  & 0  \\
   0  &  0
\end{array}\right)$.
Let $M$ be a compact manifold with boundary $\partial M$. We assume that the metric $g^{M}$ on $M$ has
the following form near the boundary
 \begin{equation}
 g^{M}=\frac{1}{h(x_{n})}g^{\partial M}+dx _{n}^{2} ,
\end{equation}
where $g^{\partial M}$ is the metric on $\partial M$. Let $U\subset
M$ be a collar neighborhood of $\partial M$ which is diffeomorphic $\partial M\times [0,1)$. By the definition of $h(x_n)\in C^{\infty}([0,1))$
and $h(x_n)>0$, there exists $\tilde{h}\in C^{\infty}\big((-\varepsilon,1)\big)$ such that $\tilde{h}|_{[0,1)}=h$ and $\tilde{h}>0$ for some
sufficiently small $\varepsilon>0$. Then there exists a metric $\hat{g}$ on $\hat{M}=M\bigcup_{\partial M}\partial M\times
(-\varepsilon,0]$ which has the form on $U\bigcup_{\partial M}\partial M\times (-\varepsilon,0 ]$
 \begin{equation}
\hat{g}=\frac{1}{\tilde{h}(x_{n})}g^{\partial M}+dx _{n}^{2} ,
\end{equation}
such that $\hat{g}|_{M}=g$.
We fix a metric $\hat{g}$ on the $\hat{M}$ such that $\hat{g}|_{M}=g$.
Note $\widetilde{D}_{F}$ is the  twisted Dirac operator on the spinor bundle $S(TM)\otimes F$ corresponding to the
connection $\widetilde{\nabla}$.

Now we recall the main theorem in \cite{FGLS}.

\begin{thm}\label{th:32}{\bf(Fedosov-Golse-Leichtnam-Schrohe)}
 Let $X$ and $\partial X$ be connected, ${\rm dim}X=n\geq3$,
 $A=\left(\begin{array}{lcr}\pi^+P+G &   K \\
T &  S    \end{array}\right)$ $\in \mathcal{B}$ , and denote by $p$, $b$ and $s$ the local symbols of $P,G$ and $S$ respectively.
 Define:
 \begin{eqnarray}
{\rm{\widetilde{Wres}}}(A)&=&\int_X\int_{\bf S}{\rm{tr}}_E\left[p_{-n}(x,\xi)\right]\sigma(\xi)dx \nonumber\\
&&+2\pi\int_ {\partial X}\int_{\bf S'}\left\{{\rm tr}_E\left[({\rm{tr}}b_{-n})(x',\xi')\right]+{\rm{tr}}
_F\left[s_{1-n}(x',\xi')\right]\right\}\sigma(\xi')dx',
\end{eqnarray}
Then

~~ a) ${\rm \widetilde{Wres}}([A,B])=0 $, for any $A,B\in\mathcal{B}$;

~~ b) It is a unique continuous trace on
$\mathcal{B}/\mathcal{B}^{-\infty}$.
\end{thm}

\section{Conformal perturbations of twisted Dirac operator and Noncommutative residue}

In this section we consider a $n$-dimensional oriented Riemannian manifold $(M, g^{M})$ equipped
with a fixed spin structure. 
Let $S(TM)$ be the spinors bundle and $F$ be an additional smooth vector bundle  equipped with a non-unitary connection $\widetilde{\nabla}^{F}$.
Let $S_1, S_2\in \Gamma(F)$, $g^F$ be a metric on $F$. We define the dual connection $\widetilde{\nabla}^{F,\ast}$ by $$g^F(\widetilde{\nabla}^{F}_XS_1,S_2)+g^F(S_1,\widetilde{\nabla}^{F,\ast}_XS_2)=X(g^F(S_1,S_2))$$
for $X\in\Gamma(TM)$ and define
 \begin{equation}
\nabla^{F}=\frac{\widetilde{\nabla}^{F}+\widetilde{\nabla}^{F,\ast}}{2},~~A=\frac{\widetilde{\nabla}^{F}-\widetilde{\nabla}^{F,\ast}}{2},
\end{equation}
then $\nabla^{F}$ is a metric connection and $\Phi$ is an endomorphism of $F$ with a 1-form coefficient.
 We consider the tensor product vector bundle $S(TM)\otimes F$,
which becomes a Clifford module via the definition:
\begin{equation}
c(a)=c(a)\otimes\rm {id}_{F},~~~~~a\in TM,
\end{equation}
and which we equip with the compound connection:
 \begin{equation}
\widetilde{\nabla}^{ S(TM)\otimes F}= \nabla^{ S(TM)}\otimes \rm {id}_{ F}+ \rm {id}_{ S(TM)}\otimes \widetilde{\nabla}^{F}.
\end{equation}

Let
 \begin{equation}
\nabla^{ S(TM)\otimes F}=\nabla^{ S(TM)}\otimes \rm {id}_{ F}+ \rm {id}_{ S(TM)}\otimes \nabla^{F},
\end{equation}
then the spinor connection  $\widetilde{\nabla}$  induced by $\nabla^{ S(TM)\otimes F}$ is locally
given by
 \begin{equation}
\widetilde{\nabla}^{ S(TM)\otimes F}=\nabla^{ S(TM)}\otimes \rm {id}_{ F}+ \rm {id}_{ S(TM)}\otimes \nabla^{F}+\rm {id}_{ S(TM)}\otimes A.
\end{equation}

Let  $\{e_{i}\}(1\leq i,j\leq n)$ $(\{\partial_{i}\})$ be the orthonormal frames (natural frames respectively ) on  $TM$,
 \begin{equation}
D_{F}=\sum^{n}_{i,j=1}g^{ij}c(\partial_{i})\nabla^{ S(TM)\otimes F}_{\partial_{j}}=\sum_{j=1}^{n}c(e_{j})\nabla^{ S(TM)\otimes F}_{e_{j}},
\end{equation}
where $\nabla^{S(TM)\otimes F}_{\partial_{j}}=\partial_{j}+\sigma_{j}^{s}+\sigma_{j}^{F}$ and
$\sigma_{j}^{s}=\frac{1}{4}\sum\limits^{n}_{j,k=1}\langle \nabla^{TM}_{\partial_{i}}e_{j}, e_{k}\rangle c(e_{j})c(e_{k})$,
$\sigma_{j}^{F}$ is the connection  matrix of $\nabla^{F}$, then the  twisted Dirac operators $\widetilde{D}_{F}$, $\widetilde{D}^{*}_{F}$ associated to the  connection
$\widetilde{\nabla}$ as follows.

For $\psi\otimes \chi\in S(TM)\otimes F$, we have
\begin{eqnarray}
&&\widetilde{D}_{F}(\psi\otimes \chi)=D_{F}(\psi\otimes \chi)+c(A)(\psi\otimes \chi),\\
&&\widetilde{D}^{*}_{F}(\psi\otimes \chi)=D_{F}(\psi\otimes \chi)-c(A^{*})(\psi\otimes \chi),
\end{eqnarray}
where $c(A)=\sum\limits_{i=1}^{n}c(e_{i})\otimes A(e_{i})$ and $c(A^{*})=\sum\limits_{i=1}^{n}c(e_{i})\otimes A^{*}(e_{i})$, $A^{*}(e_{i})$ denotes the adjoint of $A(e_{i})$.

Then, we have obtain
\begin{eqnarray}
&&\widetilde{D}_{F}=\sum_{j=1}^{n}c(e_{j})\nabla^{ S(TM)\otimes F}_{e_{j}}+c(A),\\
&&\widetilde{D}^{*}_{F}=\sum_{j=1}^{n}c(e_{j})\nabla^{ S(TM)\otimes F}_{e_{j}}-c(A^{*}).
\end{eqnarray}
Let $\nabla^{TM}$ denote the Levi-civita connection about $g^M$. In the local coordinates $\{x_i; 1\leq i\leq n\}$ and the fixed orthonormal frame
$\{\widetilde{e_1},\cdots,\widetilde{e_n}\}$, the connection matrix $(\omega_{s,t})$ is defined by
\begin{equation}
\nabla^{TM}(\widetilde{e_1},\cdots,\widetilde{e_n})= (\widetilde{e_1},\cdots,\widetilde{e_n})(\omega_{s,t}).
\end{equation}
Let $c(\widetilde{e_i})$ denote the Clifford action, $g^{ij}=g(dx_i,dx_j)$,$\nabla^{TM}_{\partial_i}\partial_j=\sum\limits_k\Gamma_{ij}^k\partial_k, ~\Gamma^k=g^{ij}\Gamma_{ij}^k$
and the cotangent vector $\xi=\sum \xi_jdx_j$ and $\xi^j=g^{ij}\xi_i$,
 by Lemma 1 in \cite{Wa4} and Lemma 2.1 in \cite{Wa3}, for any fixed point $x_0\in\partial M$, we can choose the normal coordinates $U$
 of $x_0$ in $\partial M$ (not in $M$), by the composition formula and (2.2.11) in \cite{Wa3}, we obtain in \cite{WJ6},

\begin{lem}
Let $\widetilde{D}^{*}_{F},  \widetilde{D}_{F}$ be the twisted  Dirac operators on  $\Gamma(S(TM)\otimes F)$, then
\begin{eqnarray}
&&\sigma_{-1}(\widetilde{D}^{*}_{F})^{-1}=\sigma_{-1}(\widetilde{D}^{-1}_F)=\frac{\sqrt{-1}c(\xi)}{|\xi|^2}; \\
&&\sigma_{-2}(\widetilde{D}^{*}_{F})^{-1}=\frac{c(\xi)
\sigma_{0}(\widetilde{D}^{*}_{F})c(\xi)}{|\xi|^4}+\frac{c(\xi)}{|\xi|^6}\sum_jc(dx_j)
\Big[\partial_{x_j}[c(\xi)]|\xi|^2-c(\xi)\partial_{x_j}(|\xi|^2)\Big] ;\\
&& \sigma_{-2}(\widetilde{D}^{-1}_{F})=\frac{c(\xi)\sigma_{0}(\widetilde{D}_{F})c(\xi)}{|\xi|^4}+\frac{c(\xi)}{|\xi|^6}\sum_jc(dx_j)
\Big[\partial_{x_j}[c(\xi)]|\xi|^2-c(\xi)\partial_{x_j}(|\xi|^2)\Big],
\end{eqnarray}
where
\begin{eqnarray}
\sigma_0(\widetilde{D}^{*}_{F})&=& -\frac{1}{4}\sum_{s,t}\omega_{s,t}(e_l)c(e_{l})c(e_s)c(e_t)+\sum_{j=1}^{n}c(e_{j})\big(\sigma_{j}^{F}-A^{*}(e_{j})\big);\\
\sigma_0(\widetilde{D}_{F})&=& -\frac{1}{4}\sum_{s,t}\omega_{s,t}(e_l)c(e_{l})c(e_s)c(e_t)+\sum_{j=1}^{n}c(e_{j})\big(\sigma_{j}^{F}+A(e_{j})\big).
\end{eqnarray}
\end{lem}
For convenience, let
$\lambda=\sum_{j=1}^{n}c(e_{j})\big(\sigma_{j}^{F}-A^{*}(e_{j})\big),
\mu=\sum_{j=1}^{n}c(e_{j})\big(\sigma_{j}^{F}+A(e_{j})\big)$.
In the following, we will compute the more general case $\widetilde{Wres}[\pi^{+}(f\widetilde{D}_{F}^{-1}) \circ\pi^{+}\big(f^{-1}(\widetilde{D}_{F}^{*})^{-1}\cdot f\widetilde{D}_{F}^{-1}\cdot f^{-1}(\widetilde{D}^{*}_{F})^{-1}\big)]$ for nonzero
smooth functions $f,~f^{-1}$.
Denote by $\sigma_{l}(P)$ the $l$-order symbol of an operator P.
An application of (3.5) and (3.6) in \cite{Wa5} shows that
\begin{eqnarray}
&&\widetilde{Wres}[\pi^{+}(f\widetilde{D}_{F}^{-1}) \circ\pi^{+}\big(f^{-1}(\widetilde{D}_{F}^{*})^{-1}\cdot f\widetilde{D}_{F}^{-1}\cdot f^{-1}(\widetilde{D}^{*}_{F})^{-1}\big)]\nonumber\\
&=&\int_{M}\int_{|\xi|=1}{{\rm trace}}_{S(TM)\otimes F}\big[\sigma_{-n}\big( (\widetilde{D}_{F}^{*}f\cdot\widetilde{D}_{F}f^{-1})^{-2}\big)\big]\sigma(\xi)dx+\int_{\partial M}\Phi,
\end{eqnarray}
where
 \begin{eqnarray}
\Phi&=&\int_{|\xi'|=1}\int_{-\infty}^{+\infty}\sum_{j,k=0}^{\infty}\sum \frac{(-i)^{|\alpha|+j+k+\ell}}{\alpha!(j+k+1)!}
{{\rm trace}}_{S(TM)\otimes F}\Big[\partial_{x_{n}}^{j}\partial_{\xi'}^{\alpha}\partial_{\xi_{n}}^{k}\sigma_{r}^{+}
(f\widetilde{D}_{F}^{-1})(x',0,\xi',\xi_{n})\nonumber\\
&&\times\partial_{x_{n}}^{\alpha}\partial_{\xi_{n}}^{j+1}\partial_{x_{n}}^{k}\sigma_{l}
\Big(f^{-1}(\widetilde{D}_{F}^{*})^{-1}\cdot f\widetilde{D}_{F}^{-1}\cdot f^{-1}(\widetilde{D}^{*}_{F})^{-1}\Big)(x',0,\xi',\xi_{n})\Big]
d\xi_{n}\sigma(\xi')\texttt{d}x' ,
\end{eqnarray}
and the sum is taken over $r-k+|\alpha|+\ell-j-1=-n=-6,r\leq-1,\ell\leq-3$.

Note that
\begin{eqnarray}
&&f^{-1}(\widetilde{D}_{F}^{*})^{-1}\cdot f\widetilde{D}_{F}^{-1}\cdot f^{-1}(\widetilde{D}^{*}_{F})^{-1}\nonumber\\
&=&(\widetilde{D}^{*}_{F}f\cdot
\widetilde{D}_{F}f^{-1}\cdot\widetilde{D}_{F}^{*}f)^{-1}\nonumber\\
&=&\Big(\widetilde{D}^{*}_{F}f\cdot
\widetilde{D}_{F}\widetilde{D}_{F}^{*}f^{-1}\cdot f-\widetilde{D}^{*}_{F}f\cdot
\widetilde{D}_{F}\cdot[\widetilde{D}_{F}^{*},f^{-1}]\cdot f\Big)^{-1}\nonumber\\
&=&\Big(\widetilde{D}^{*}_{F}f\cdot
\widetilde{D}_{F}\widetilde{D}_{F}^{*}-\widetilde{D}^{*}_{F}f\cdot
\widetilde{D}_{F}c(df^{-1})f\Big)^{-1}\nonumber\\
&=&\Big(f\cdot\widetilde{D}_{F}^{*}\widetilde{D}_{F}\widetilde{D}_{F}^{*}
+[\widetilde{D}_{F}^{*},f]\widetilde{D}_{F}\widetilde{D}_{F}^{*}
-\widetilde{D}_{F}^{*}f\cdot\widetilde{D}_{F}c(df^{-1})f\Big)^{-1}\nonumber\\
&=&\Big(f\cdot\widetilde{D}_{F}^{*}\widetilde{D}_{F}\widetilde{D}_{F}^{*}
+c(df)\widetilde{D}_{F}\widetilde{D}_{F}^{*}
-\widetilde{D}_{F}^{*}f\cdot\widetilde{D}_{F}c(df^{-1})f\Big)^{-1}\nonumber\\
&=&\Big(f\cdot\widetilde{D}_{F}^{*}\widetilde{D}_{F}\widetilde{D}_{F}^{*}
+c(df)\widetilde{D}_{F}\widetilde{D}_{F}^{*}
-\widetilde{D}_{F}^{*}\widetilde{D}_{F}f\cdot c(df^{-1})\cdot f
+\widetilde{D}_{F}^{*}\cdot c(df) c(df^{-1})f\Big)^{-1}.
\end{eqnarray}
In order to get the symbol of operators $\widetilde{D}^{*}_{F}f\cdot
\widetilde{D}_{F}f^{-1}\cdot\widetilde{D}_{F}^{*}f$. We first give the specification of
$\widetilde{D}^{*}_{F}\widetilde{D}_{F}\widetilde{D}_{F}^{*}$, $\widetilde{D}^{*}_{F}\widetilde{D}_{F}$ and $\widetilde{D}_{F}\widetilde{D}_{F}^{*}$.
By (3.9) and (3.10), we have
\begin{eqnarray}
&&\widetilde{D}_{F}\widetilde{D}^{*}_{F}\nonumber\\
&=&D_{F}^{2}-D_{F}c(A^{*})+c(A)D_{F}-c(A)c(A^{*})\nonumber\\
                &=&-g^{ij}\partial_{i}\partial_{j}-2\sigma^{j}_{S(TM)\otimes F}\partial_{j}+\Gamma^{k}\partial_{k}
                +\sum^{n}_{j=1}\Big[c(A)c(e_{j})-c(e_{j})c(A^{*}) \Big]e_{j}-\sum^{n}_{j=1}c(e_{j})\sigma_{j}^{S(TM)\otimes F} c(A^*)\nonumber\\
                &&-g^{ij}\Big[(\partial_{i}\sigma^{j}_{S(TM)\otimes F}) +\sigma^{i}_{S(TM)\otimes F}\sigma^{j}_{S(TM)\otimes F}
                  -\Gamma_{ij}^{k}\sigma_{S(TM)\otimes F}^{k}\Big]+\frac{1}{4}s+\frac{1}{2}\sum_{i\neq j} R^{F}(e_{i},e_{j})c(e_{i})c(e_{j})\nonumber\\
                &&+\sum^{n}_{j=1}\Big[c(A)c(e_{j}) \Big]\sigma^{S(TM)\otimes F}_{j}
               -\sum^{n}_{j=1}c(e_{j}) e_{j}\big(c(A^{*})\big)-c(A)c(A^{*})
\end{eqnarray}
and
\begin{eqnarray}
&&\widetilde{D}^{*}_{F}\widetilde{D}_{F}\nonumber\\
&=&D_{F}^{2}-c(A^{*})D_{F}+D_{F}c(A)-c(A)c(A^{*})\nonumber\\
                &=&-g^{ij}\partial_{i}\partial_{j}-2\sigma^{j}_{S(TM)\otimes F}\partial_{j}+\Gamma^{k}\partial_{k}
                +\sum^{n}_{j=1}\Big[c(e_{j})c(A)-c(A^{*}) c(e_{j}) \Big]e_{j}+\sum^{n}_{j=1}c(e_{j})\sigma_{j}^{S(TM)\otimes F} c(A)\nonumber\\
                &&-g^{ij}\Big[(\partial_{i}\sigma^{j}_{S(TM)\otimes F}) +\sigma^{i}_{S(TM)\otimes F}\sigma^{j}_{S(TM)\otimes F}
                  -\Gamma_{ij}^{k}\sigma_{S(TM)\otimes F}^{k}\Big]+\frac{1}{4}s+\frac{1}{2}\sum_{i\neq j} R^{F}(e_{i},e_{j})c(e_{i})c(e_{j})\nonumber\\
                &&-\sum^{n}_{j=1}\Big[c(A^{*})c(e_{j}) \Big]\sigma^{S(TM)\otimes F}_{j}
               +\sum^{n}_{j=1}c(e_{j}) e_{j}\big(c(A)\big)-c(A^{*})c(A).
\end{eqnarray}
Combining (3.10) and (3.20), we obtain
\begin{eqnarray}
&&\widetilde{D}^{*}_{F}\widetilde{D}_{F}\widetilde{D}^{*}_{F}\nonumber\\
         &=&\sum^{n}_{i,j,l=1}\sum^{n}_{r=1}c(e_{r})\langle e_{r},dx_{l}\rangle(-g^{ij}\partial_{l}\partial_{i}\partial_{j})
         +\sum^{n}_{r,l=1}c(e_{r})\langle e_{r},dx_{l}\rangle \bigg\{-\sum^{n}_{i,j=1}(\partial_{l}g^{ij})\partial_{i}\partial_{j}-\sum^{n}_{i,j,k=1}g^{ij}\nonumber\\
         &&\times(4\sigma^{S(TM)\otimes F}_{i}\partial_{j}-2\Gamma^{k}_{ij}\partial_{k})\partial_{l}\bigg\}
         +\sum^{n}_{r,l=1}c(e_{r})\langle e_{r},dx_{l}\rangle \bigg\{-2\sum^{n}_{i,j=1}(\partial_{l}g^{ij})\sigma^{S(TM)\otimes F}_{i}\partial_{j}\nonumber\\
         &&+\sum^{n}_{i,j,k=1}g^{ij}(\partial_{l}\Gamma^{k}_{ij})\partial_{k}
         -2\sum^{n}_{i,j=1}g^{ij}(\partial_{l}\sigma^{S(TM)\otimes F}_{i})\partial_{j}  +\sum^{n}_{i,j,k=1}(\partial_{l}g^{ij})\Gamma^{k}_{ij}\partial_{k}
         +\sum^{n}_{j,k=1}\Big[\partial_{l}\Big(c(A)c(e_{j})\nonumber\\
         &&-c(e_{j})c(A^{*})\Big)\Big]\langle e_{j},dx^{k}\rangle\partial_{k}
         +\sum^{n}_{j,k=1}\Big(c(A)c(e_{j})-c(e_{j}) c(A^{*})\Big)\Big[\partial_{l}\langle e_{j},dx^{k}\rangle\Big]\partial_{k} \bigg\}+\sum^{n}_{r,l=1}c(e_{r})\langle e_{r},dx_{l}\rangle\nonumber\\
         &&\times\partial_{l}\bigg\{-\sum^{n}_{i,j,k=1}g^{ij}\Big[(\partial_{i}\sigma^{j}_{S(TM)\otimes F}) +\sigma^{i}_{S(TM)\otimes F}\sigma^{j}_{S(TM)\otimes F}-\Gamma_{ij}^{k}\sigma_{S(TM)\otimes F}^{k}\Big]+\frac{1}{4}s+\sum^{n}_{j=1}\Big[c(A)c(e_{j}) \Big]\nonumber\\
         &&\times\sigma^{S(TM)\otimes F}_{j}-\sum^{n}_{j=1}c(e_{j}) e_{j}\big(c(A^{*})\big)
          -\sum^{n}_{j=1}c(e_{j})\sigma_{j}^{S(TM)\otimes F} c(A^{*})-c(A)c(A^{*})+\frac{1}{2}\sum_{i\neq j} R^{F}(e_{i},e_{j})\nonumber\\
         &&\times c(e_{i})c(e_{j})\bigg\}+\sigma_{0}(\widetilde{D}^{*}_{F})
         \sum^{n}_{j,i=1}(-g^{ij}\partial_{i}\partial_{j})+\sum^{n}_{r,l=1}c(e_{r})\langle e_{r},dx_{l}\rangle \bigg\{2\sum^{n}_{j,k=1}\Big[c(A)c(e_{j})-c(e_{j})c(A^{*}) \nonumber\\
         &&\Big]\langle e_{i},dx_{k}\rangle\bigg\}\partial_{l}\partial_{k}
           +\sigma_{0}(\widetilde{D}^{*}_{F})\bigg\{-2\sigma^{j}_{S(TM)\otimes F}\partial_{j}+\Gamma^{k}\partial_{k}
         +\sum^{n}_{j=1}\Big[c(A)c(e_{j})-c(e_{j}) c(A^{*}) \Big]e_{j}-\sum^{n}_{j=1}c(e_{j})\nonumber\\
         && e_{j}\big(c(A^{*})\big)
         -g^{ij}\Big[(\partial_{i}\sigma^{j}_{S(TM)\otimes F}) +\sigma^{i}_{S(TM)\otimes F}\sigma^{j}_{S(TM)\otimes F}
         -\Gamma_{ij}^{k}\sigma_{S(TM)\otimes F}^{k}\Big]+\frac{1}{4}s-c(A) c(A^{*})\nonumber\\
         &&+\sum^{n}_{j=1}\Big[c(A)c(e_{j}) \Big]\sigma^{S(TM)\otimes F}_{j}
         -\sum^{n}_{j=1}c(e_{j})\sigma_{j}^{S(TM)\otimes F} c(A^{*})+\frac{1}{2}\sum_{i\neq j} R^{F}(e_{i},e_{j})c(e_{i})c(e_{j})\bigg\}.
\end{eqnarray}

Thus, using (3.19)-(3.22), we get the specification of
$\widetilde{D}^{*}_{F}f\cdot
\widetilde{D}_{F}f^{-1}\cdot\widetilde{D}_{F}^{*}f$.

\begin{eqnarray}
&&\widetilde{D}^{*}_{F}f\cdot
\widetilde{D}_{F}f^{-1}\cdot\widetilde{D}_{F}^{*}f\nonumber\\
&=&f\cdot\widetilde{D}_{F}^{*}\widetilde{D}_{F}\widetilde{D}_{F}^{*}
+c(df)\widetilde{D}_{F}\widetilde{D}_{F}^{*}
-\widetilde{D}_{F}^{*}\widetilde{D}_{F}f\cdot c(df^{-1})\cdot f
+\widetilde{D}_{F}^{*}\cdot c(df) c(df^{-1})f\nonumber\\
&=&f\cdot\Bigg\{\sum^{n}_{i,j,l=1}\sum^{n}_{r=1}c(e_{r})\langle e_{r},dx_{l}\rangle(-g^{ij}\partial_{l}\partial_{i}\partial_{j})
         +\sum^{n}_{r,l=1}c(e_{r})\langle e_{r},dx_{l}\rangle \bigg\{-\sum^{n}_{i,j=1}(\partial_{l}g^{ij})\partial_{i}\partial_{j}-\sum^{n}_{i,j,k=1}g^{ij}\nonumber\\
         &&\times(4\sigma^{S(TM)\otimes F}_{i}\partial_{j}-2\Gamma^{k}_{ij}\partial_{k})\partial_{l}\bigg\}
         +\sum^{n}_{r,l=1}c(e_{r})\langle e_{r},dx_{l}\rangle \bigg\{-2\sum^{n}_{i,j=1}(\partial_{l}g^{ij})\sigma^{S(TM)\otimes F}_{i}\partial_{j}+\sum^{n}_{i,j,k=1}g^{ij}\nonumber\\
         &&\times(\partial_{l}\Gamma^{k}_{ij})\partial_{k}
         -2\sum^{n}_{i,j=1}g^{ij}(\partial_{l}\sigma^{S(TM)\otimes F}_{i})\partial_{j}  +\sum^{n}_{i,j,k=1}(\partial_{l}g^{ij})\Gamma^{k}_{ij}\partial_{k}
         +\sum^{n}_{j,k=1}\Big[\partial_{l}\Big(c(A)c(e_{j})-c(e_{j})\nonumber\\
         &&\times c(A^{*})\Big)\Big]\langle e_{j},dx^{k}\rangle\partial_{k}
         +\sum^{n}_{j,k=1}\Big(c(A)c(e_{j})-c(e_{j}) c(A^{*})\Big)\Big[\partial_{l}\langle e_{j},dx^{k}\rangle\Big]\partial_{k} \bigg\}+\sum^{n}_{r,l=1}c(e_{r})\langle e_{r},dx_{l}\rangle\nonumber\\
         &&\times\partial_{l}\bigg\{-\sum^{n}_{i,j,k=1}g^{ij}\Big[(\partial_{i}\sigma^{j}_{S(TM)\otimes F}) +\sigma^{i}_{S(TM)\otimes F}\sigma^{j}_{S(TM)\otimes F}-\Gamma_{ij}^{k}\sigma_{S(TM)\otimes F}^{k}\Big]+\frac{1}{4}s+\sum^{n}_{j=1}\Big[c(A)\nonumber\\
         &&\times c(e_{j}) \Big]\sigma^{S(TM)\otimes F}_{j}-\sum^{n}_{j=1}c(e_{j}) e_{j}\big(c(A^{*})\big)
          -\sum^{n}_{j=1}c(e_{j})\sigma_{j}^{S(TM)\otimes F} c(A^{*})-c(A)c(A^{*})+\frac{1}{2}\sum_{i\neq j} R^{F}(e_{i},e_{j})\nonumber\\
         &&\times c(e_{i})c(e_{j})\bigg\}+\sigma_{0}(\widetilde{D}^{*}_{F})
         \sum^{n}_{j,i=1}(-g^{ij}\partial_{i}\partial_{j})+\sum^{n}_{r,l=1}c(e_{r})\langle e_{r},dx_{l}\rangle \bigg\{2\sum^{n}_{j,k=1}\Big[c(A)c(e_{j})-c(e_{j})c(A^{*})\Big] \nonumber\\
         &&\times\langle e_{i},dx_{k}\rangle\bigg\}\partial_{l}\partial_{k}
           +\sigma_{0}(\widetilde{D}^{*}_{F})\bigg\{-2\sigma^{j}_{S(TM)\otimes F}\partial_{j}+\Gamma^{k}\partial_{k}
         +\sum^{n}_{j=1}\Big[c(A)c(e_{j})-c(e_{j}) c(A^{*}) \Big]e_{j}
         -\sum^{n}_{j=1}c(e_{j})\nonumber\\
         &&\times e_{j}\big(c(A^{*})\big)
         -g^{ij}\Big[(\partial_{i}\sigma^{j}_{S(TM)\otimes F}) +\sigma^{i}_{S(TM)\otimes F}\sigma^{j}_{S(TM)\otimes F}
         -\Gamma_{ij}^{k}\sigma_{S(TM)\otimes F}^{k}\Big]+\frac{1}{4}s-c(A)c(A^{*})\nonumber\\
         &&+\sum^{n}_{j=1}\Big[c(A)c(e_{j}) \Big]\sigma^{S(TM)\otimes F}_{j}
         -\sum^{n}_{j=1}c(e_{j})\sigma_{j}^{S(TM)\otimes F} c(A^{*})+\frac{1}{2}\sum_{i\neq j} R^{F}(e_{i},e_{j})c(e_{i})c(e_{j})\bigg\}\Bigg\}+c(df)\cdot\nonumber\\
         &&\Bigg\{
         -g^{ij}\partial_{i}\partial_{j}-2\sigma^{j}_{S(TM)\otimes F}\partial_{j}+\Gamma^{k}\partial_{k}
                +\sum^{n}_{j=1}\Big[c(A)c(e_{j})-c(e_{j})c(A^{*}) \Big]e_{j}-\sum^{n}_{j=1}c(e_{j})\sigma_{j}^{S(TM)\otimes F}\nonumber\\
                &&\times c(A^*)-g^{ij}\Big[(\partial_{i}\sigma^{j}_{S(TM)\otimes F}) +\sigma^{i}_{S(TM)\otimes F}\sigma^{j}_{S(TM)\otimes F}
                  -\Gamma_{ij}^{k}\sigma_{S(TM)\otimes F}^{k}\Big]+\frac{1}{4}s+\frac{1}{2}\sum_{i\neq j} R^{F}(e_{i},e_{j})\nonumber\\
                &&\times c(e_{i})c(e_{j})+\sum^{n}_{j=1}\Big[c(A)c(e_{j}) \Big]\sigma^{S(TM)\otimes F}_{j}
               -\sum^{n}_{j=1}c(e_{j}) e_{j}\big(c(A^{*})\big)-c(A)c(A^{*})\Bigg\}-
\Bigg\{-g^{ij}\partial_{i}\nonumber\\
               &&\times\partial_{j}-2\sigma^{j}_{S(TM)\otimes F}\partial_{j}+\Gamma^{k}\partial_{k}
                +\sum^{n}_{j=1}\Big[c(e_{j})c(A)-c(A^{*}) c(e_{j}) \Big]e_{j}+\sum^{n}_{j=1}c(e_{j})\sigma_{j}^{S(TM)\otimes F} c(A)\nonumber\\
                &&-g^{ij}\Big[(\partial_{i}\sigma^{j}_{S(TM)\otimes F}) +\sigma^{i}_{S(TM)\otimes F}\sigma^{j}_{S(TM)\otimes F}
                  -\Gamma_{ij}^{k}\sigma_{S(TM)\otimes F}^{k}\Big]+\frac{1}{4}s+\frac{1}{2}\sum_{i\neq j} R^{F}(e_{i},e_{j})c(e_{i})\nonumber\\
                &&\times c(e_{j})-\sum^{n}_{j=1}\Big[c(A^{*})c(e_{j}) \Big]\sigma^{S(TM)\otimes F}_{j}
               +\sum^{n}_{j=1}c(e_{j}) e_{j}\big(c(A)\big)-c(A^{*})c(A)\Bigg\}f\cdot c(df^{-1})\cdot f\nonumber\\
               &&+\Bigg\{\sum^{n}_{i,j=1}g^{ij}c(\partial_{i})\Big(\partial_{j}+\sigma^{S(TM)\otimes F}_{j}\Big)-c(A^*)\Bigg\}\cdot c(df) c(df^{-1})f.
    \end{eqnarray}

Let $\partial^{j}=g^{ij}\partial_{i}, \sigma^{i}=g^{ij}\sigma_{j}$, 
by the above formulas, then we obtain:

\begin{lem}
Let $\widetilde{D}^{*}_{F}, \widetilde{D}_{F}$ be the twisted  Dirac operators on  $\Gamma(S(TM)\otimes F)$,
\begin{eqnarray}
\sigma_{3}(\widetilde{D}^{*}_{F}f\cdot\widetilde{D}_{F}f^{-1}\cdot\widetilde{D}^{*}_{F}f)
&=& f\sigma_{3}(\widetilde{D}^{*}_{F}\widetilde{D}_{F}\widetilde{D}^{*}_{F})
=\sqrt{-1}c(\xi)|\xi|^2f; \\
\sigma_{2}(\widetilde{D}^{*}_{F}f\cdot\widetilde{D}_{F}f^{-1}\cdot\widetilde{D}^{*}_{F}f)
&=&f\sigma_{2}(\widetilde{D}^{*}_{F}\widetilde{D}_{F}\widetilde{D}^{*}_{F})
+2c(df)|\xi|^2,
\end{eqnarray}
where
$
\sigma_{2}(\widetilde{D}^{*}_{F}\widetilde{D}_{F}\widetilde{D}^{*}_{F})
=c(\xi)(4\sigma^k-2\Gamma^k)\xi_{k}-\frac{1}{4}|\xi|^2h'(0)c(dx_{n})
+\lambda|\xi|^2-2c(\xi)c(A)c(\xi)-2|\xi|^2c(A^*)$.
\end{lem}
For convenience,
we write that
$\sigma_{2}(\widetilde{D}^{*}_{F}\widetilde{D}_{F}\widetilde{D}^{*}_{F})=
G+\lambda|\xi|^2-2c(\xi)c(A)c(\xi)-2|\xi|^2c(A^*)$.
In order to get the symbol of operators $\widetilde{D}^{*}_{F}f\cdot
\widetilde{D}_{F}f^{-1}\cdot\widetilde{D}_{F}^{*}f$. We first give the following  formulas:
\begin{eqnarray}
D_x^{\alpha}=(-\sqrt{-1})^{|\alpha|}\partial_x^{\alpha};
~\sigma(\widetilde{D}^{*}_{F}f\cdot\widetilde{D}_{F}f^{-1}\cdot\widetilde{D}^{*}_{F}f)&=&p_3+p_2+p_1+p_0;\nonumber\\
~\sigma\big((\widetilde{D}^{*}_{F}f\cdot\widetilde{D}_{F}f^{-1}
\cdot\widetilde{D}^{*}_{F}f)^{-1}\big)&=&\sum^{\infty}_{j=3}q_{-j}.
\end{eqnarray}
By the composition formula of psudodifferential operators, we have
 \begin{eqnarray}
1
&=&\sigma\big[(\widetilde{D}^{*}_{F}f\cdot\widetilde{D}_{F}f^{-1}\cdot\widetilde{D}^{*}_{F}f)\circ (\widetilde{D}^{*}_{F}f\cdot\widetilde{D}_{F}f^{-1}\cdot\widetilde{D}^{*}_{F}f)^{-1}\big]\nonumber\\
&=&(p_3+p_2+p_1+p_0)(q_{-3}+q_{-4}+q_{-5}+\cdots) \nonumber\\
&&+\sum_j(\partial_{\xi_j}p_3+\partial_{\xi_j}p_2+\partial_{\xi_j}p_1+\partial_{\xi_j}p_0)
(D_{x_j}q_{-3}+D_{x_j}q_{-4}+D_{x_j}q_{-5}+\cdots) \nonumber\\
&=&p_3q_{-3}+(p_3q_{-4}+p_2q_{-3}+\sum_j\partial_{\xi_j}p_3D_{x_j}q_{-3})+\cdots.
\end{eqnarray}
Then
\begin{equation}
q_{-3}=p_3^{-1};~q_{-4}=-p_3^{-1}\Big[p_2p_3^{-1}+\sum_j\partial_{\xi_j}p_3D_{x_j}(p_3^{-1})\Big].
\end{equation}\\
By Lemma 2.1 in \cite{Wa3} and (3.24), (3.25), we obtain

\begin{lem}Let $\widetilde{D}^{*}_{F}, \widetilde{D}_{F}$ be the twisted  Dirac operators on  $\Gamma(S(TM)\otimes F)$, then
\begin{eqnarray}
\sigma_{-3}(\widetilde{D}^{*}_{F}f\cdot\widetilde{D}_{F}f^{-1}\cdot\widetilde{D}^{*}_{F}f)^{-1}
&=&f^{-1}\sigma_{-3}(\widetilde{D}^{*}_{F}\widetilde{D}_{F}\widetilde{D}^{*}_{F})^{-1}
=\frac{\sqrt{-1}c(\xi)}{f|\xi|^4}; \\
\sigma_{-4}(\widetilde{D}^{*}_{F}f\cdot\widetilde{D}_{F}f^{-1}\cdot\widetilde{D}^{*}_{F}f)^{-1}
&=&f^{-1}\sigma_{-4}(\widetilde{D}^{*}_{F}\widetilde{D}_{F}\widetilde{D}^{*}_{F})^{-1}
+\frac{2c(\xi)c(df)c(\xi)}{f^{2}|\xi|^6}\nonumber\\
&&+\frac{ic(\xi)\sum\limits_j\Big[c(dx_j)|\xi|^2+2\xi_{j}c(\xi)\Big]D_{x_j}(f^{-1})c(\xi)}{|\xi|^8},
\end{eqnarray}
where
\begin{eqnarray}
&&\sigma_{-4}(\widetilde{D}^{*}_{F}\widetilde{D}_{F}\widetilde{D}^{*}_{F})^{-1}\nonumber\\
&=&\frac{c(\xi)\sigma_{2}(\widetilde{D}^{*}_{F}\widetilde{D}_{F}\widetilde{D}^{*}_{F})c(\xi)}{|\xi|^8}
+\frac{c(\xi)}{|\xi|^{10}}\sum_j\Big[c(dx_j)|\xi|^2
+2\xi_{j}c(\xi)\Big]\Big[\partial_{x_j}[c(\xi)]|\xi|^2-2c(\xi)\partial_{x_j}(|\xi|^2)\Big]\nonumber\\
&=&\frac{c(\xi)Gc(\xi)}{|\xi|^8}
+\frac{c(\xi)\lambda c(\xi)}{|\xi|^6}-\frac{2c(A)}{|\xi|^4}-\frac{2c(\xi)c(A^*)c(\xi)}{|\xi|^6}
+\frac{c(\xi)}{|\xi|^{10}}\sum_j\Big[c(dx_j)|\xi|^2
+2\xi_{j}c(\xi)\Big]\Big[\partial_{x_j}[c(\xi)]|\xi|^2\nonumber\\
&&-2c(\xi)\partial_{x_j}(|\xi|^2)\Big]
\end{eqnarray}
\end{lem}

Locally we can use Theorem 2.5 in [19] to compute the interior term of (3.17), then
 \begin{eqnarray}
&&\int_{M}\int_{|\xi|=1}{{\rm trace}}_{S(TM)\otimes F}[\sigma_{-n}\big( (\widetilde{D}_{F}^{*}f\cdot\widetilde{D}_{F}f^{-1})^{-2}\big)\big]\sigma(\xi)dx\nonumber\\
&=& 8 \pi^{3}\int_{M}\bigg \{{\rm{trace}}\Big[-\frac{s}{12}+c(A^{*})c(A)-\frac{1}{4}\sum_{i}\big[c(A^{*})c(e_{i})
-c(e_{i})c(A) \big]^{2}\nonumber\\
 &&-\frac{1}{2}\sum_{j}\nabla_{e_{j}}^{F}\big(c(A^{*})\big)c(e_{j})
 -\frac{1}{2}\sum_{j}c(e_{j})\nabla_{e_{j}}^{F}\big(c(A)\big)\Big]
 -2f^{-1}\Delta(f)\nonumber\\
 &&+4f^{-1}{\rm{trace}}\Big[A(grad_{M}f)\Big]-f^{2}\Big[|grad_{M}(f)|^{2}+2\Delta(f)\Big]\bigg \}dvol_{M}.
\end{eqnarray}
So we only need to compute $\int_{\partial M}\Phi$.

 From the formula (3.18) for the definition of $\Phi$, now we can compute $\Phi$.
Since the sum is taken over $r+\ell-k-j-|\alpha|-1=-6, \ r\leq-1, \ell\leq -3$, then we have the $\int_{\partial_{M}}\Phi$
is the sum of the following five cases:
~\\
~\\
\noindent  {\bf case (a)~(I)}~$r=-1, l=-3, j=k=0, |\alpha|=1$.\\

By (3.18), we get
 \begin{eqnarray}
{\rm case~(a)~(I)}&=&-\int_{|\xi'|=1}\int^{+\infty}_{-\infty}\sum_{|\alpha|=1}{\rm trace}
\Big[\partial^{\alpha}_{\xi'}\pi^{+}_{\xi_{n}}\sigma_{-1}(f\widetilde{D}_{F}^{-1})
      \times\partial^{\alpha}_{x'}\partial_{\xi_{n}}\sigma_{-3}
      \big(f^{-1}(\widetilde{D}_{F}^{*})^{-1}\cdot f\widetilde{D}_{F}^{-1}\cdot\nonumber\\
      &&f^{-1}(\widetilde{D}^{*}_{F})^{-1}\big)\Big](x_0)d\xi_n\sigma(\xi')dx'\nonumber\\
&=&-\int_{|\xi'|=1}\int^{+\infty}_{-\infty}\sum_{|\alpha|=1}{\rm trace}
\Big[\partial^{\alpha}_{\xi'}\pi^{+}_{\xi_{n}}\big(f\sigma_{-1}(\widetilde{D}_{F}^{-1})\big)
      \times\partial^{\alpha}_{x'}\partial_{\xi_{n}}\big(f^{-1}\sigma_{-3}
(\widetilde{D}_{F}^{*}\widetilde{D}_{F}\widetilde{D}^{*}_{F})^{-1}\big)\Big](x_0)\nonumber\\
      &&\times d\xi_n\sigma(\xi')dx'\nonumber\\
  &=&-\int_{|\xi'|=1}\int^{+\infty}_{-\infty}\sum_{|\alpha|=1}{\rm trace}
\Big[\partial^{\alpha}_{\xi'}\pi^{+}_{\xi_{n}}\sigma_{-1}(\widetilde{D}_{F}^{-1})
      \times\partial^{\alpha}_{x'}\partial_{\xi_{n}}\sigma_{-3}
(\widetilde{D}_{F}^{*}\widetilde{D}_{F}\widetilde{D}^{*}_{F})^{-1}\Big](x_0)d\xi_n\sigma(\xi')dx'\nonumber\\
      &&-f\sum\limits_{j<n}\partial_{j}(f^{-1})\int_{|\xi'|=1}\int^{+\infty}_{-\infty}\sum_{|\alpha|=1}{\rm trace}
\Big[\partial^{\alpha}_{\xi'}\pi^{+}_{\xi_{n}}\sigma_{-1}(\widetilde{D}_{F}^{-1})
      \times\partial_{\xi_{n}}\sigma_{-3}
(\widetilde{D}_{F}^{*}\widetilde{D}_{F}\widetilde{D}^{*}_{F})^{-1}\Big](x_0)\nonumber\\
      &&\times d\xi_n\sigma(\xi')dx'.
\end{eqnarray}
By Lemma 2.2 in \cite{Wa3} and (3.29), for $i<n$, we have
 \begin{eqnarray}
 &&\partial_{x_{i}}\sigma_{-3}\Big((\widetilde{D}^{*}_{F}\widetilde{D}_{F}\widetilde{D}^{*}_{F})^{-1}\Big)(x_0)=
      \partial_{x_{i}}\Big[\frac{\sqrt{-1}c(\xi)}{|\xi|^{4}}\Big](x_{0})\nonumber\\
      &=&\sqrt{-1}\partial_{x_{i}}\Big[c(\xi)\Big]|\xi|^{-4}(x_{0})
      -2\sqrt{-1}c(\xi)\partial_{x_{i}}\Big[|\xi|^{2}\Big]|\xi|^{-6}(x_{0})=0.
\end{eqnarray}
Thus we have
\begin{eqnarray}
-\int_{|\xi'|=1}\int^{+\infty}_{-\infty}\sum_{|\alpha|=1}{\rm trace}
\Big[\partial^{\alpha}_{\xi'}\pi^{+}_{\xi_{n}}\sigma_{-1}(\widetilde{D}_{F}^{-1})
      \times\partial^{\alpha}_{x'}\partial_{\xi_{n}}\sigma_{-3}
(\widetilde{D}_{F}^{*}\widetilde{D}_{F}\widetilde{D}^{*}_{F})^{-1}\Big](x_0)d\xi_n\sigma(\xi')dx'=0.
\end{eqnarray}
By (3.12) and direct calculations, for $i<n$, we obtain
\begin{eqnarray}
&&\partial^{\alpha}_{\xi'}\pi^{+}_{\xi_{n}}\sigma_{-1}(\widetilde{D}_{F}^{-1})(x_0)|_{|\xi'|=1}
=\partial_{\xi_i}\pi^{+}_{\xi_{n}}\sigma_{-1}(\widetilde{D}_{F}^{-1})(x_0)|_{|\xi'|=1}\nonumber\\
&=&\frac{c(dx_i)}{2(\xi_n-\sqrt{-1})}-\frac{\xi_i(\xi_n-2\sqrt{-1})c(\xi')+\xi_ic(dx_n)}{2(\xi_n-\sqrt{-1})^2},
\end{eqnarray}
and we get
\begin{eqnarray}
\partial_{\xi_{n}}\sigma_{-3}(\widetilde{D}_{F}^{*}\widetilde{D}_{F}\widetilde{D}^{*}_{F})^{-1}
=\frac{\sqrt{-1}c(dx_n)}{|\xi|^4}-\frac{4\sqrt{-1}\big[\xi_nc(\xi')+\xi^2_nc(dx_n)\big]}{|\xi|^6}.
\end{eqnarray}
Then for $i<n$, we have
\begin{eqnarray}
&&{\rm trace}
\Big[\partial^{\alpha}_{\xi'}\pi^{+}_{\xi_{n}}\sigma_{-1}(\widetilde{D}_{F}^{-1})
      \times\partial_{\xi_{n}}\sigma_{-3}
(\widetilde{D}_{F}^{*}\widetilde{D}_{F}\widetilde{D}^{*}_{F})^{-1}\Big](x_0)\nonumber\\
&=&-\xi_i{\rm trace}
\Big[\frac{c(dx_n)^{2}}{2(\xi_n-\sqrt{-1})^2}\Big]-4\sqrt{-1}\xi_n\xi_i{\rm trace}
\Big[\frac{c(dx_i)^{2}}{2(\xi_n-\sqrt{-1})|\xi|^6}\Big]+4\sqrt{-1}\xi_n\xi_i(\xi_n-2\sqrt{-1})\nonumber\\
&&\times{\rm trace}
\Big[\frac{c(\xi')^{2}}{2(\xi_n-\sqrt{-1})^2|\xi|^6}\Big]+4\sqrt{-1}\xi^{2}_n\xi_i{\rm trace}
\Big[\frac{c(dx_n)^{2}}{2(\xi_n-\sqrt{-1})^2|\xi|^6}\Big].
\end{eqnarray}
We note that $i<n,~\int_{|\xi'|=1}\xi_i\sigma(\xi')=0$,
so
\begin{eqnarray}
&&-f\sum\limits_{j<n}\partial_{j}(f^{-1})\int_{|\xi'|=1}\int^{+\infty}_{-\infty}\sum_{|\alpha|=1}{\rm trace}
\Big[\partial^{\alpha}_{\xi'}\pi^{+}_{\xi_{n}}\sigma_{-1}(\widetilde{D}_{F}^{-1})
      \times\partial_{\xi_{n}}\sigma_{-3}
(\widetilde{D}_{F}^{*}\widetilde{D}_{F}\widetilde{D}^{*}_{F})^{-1}\Big](x_0) d\xi_n\sigma(\xi')dx'\nonumber\\
&=&0.
\end{eqnarray}
Then we have ${\bf case~(a)~(I)}=0$.
~\\

\noindent  {\bf case (a)~(II)}~$r=-1, l=-3, |\alpha|=k=0, j=1$.\\

By (3.18), we have
  \begin{eqnarray}
{\rm case~(a)~(II)}&=&-\frac{1}{2}\int_{|\xi'|=1}\int^{+\infty}_{-\infty} {\rm
trace} \Big[\partial_{x_{n}}\pi^{+}_{\xi_{n}}\sigma_{-1}(f\widetilde{D}_{F}^{-1})
\times\partial^{2}_{\xi_{n}}\sigma_{-3}
\big(f^{-1}(\widetilde{D}_{F}^{*})^{-1}\cdot f\widetilde{D}_{F}^{-1}\cdot\nonumber\\
&&f^{-1}(\widetilde{D}^{*}_{F})^{-1}\big)\Big](x_0)d\xi_n\sigma(\xi')dx'\nonumber\\
  &=&-\frac{1}{2}\int_{|\xi'|=1}\int^{+\infty}_{-\infty}{\rm trace}
\Big[\partial_{x_{n}}\pi^{+}_{\xi_{n}}\sigma_{-1}(\widetilde{D}_{F}^{-1})
      \times\partial^2_{\xi_{n}}\sigma_{-3}
(\widetilde{D}_{F}^{*}\widetilde{D}_{F}\widetilde{D}^{*}_{F})^{-1}\Big](x_0)d\xi_n\sigma(\xi')dx'\nonumber\\
      &&-\frac{1}{2}f^{-1}\partial_{x_{n}}(f)
      \int_{|\xi'|=1}\int^{+\infty}_{-\infty}{\rm trace}
\Big[\pi^{+}_{\xi_{n}}\sigma_{-1}(\widetilde{D}_{F}^{-1})
      \times\partial^2_{\xi_{n}}\sigma_{-3}
(\widetilde{D}_{F}^{*}\widetilde{D}_{F}\widetilde{D}^{*}_{F})^{-1}\Big](x_0)\nonumber\\
      &&\times d\xi_n\sigma(\xi')dx'.
\end{eqnarray}
By (2.2.23) in \cite{Wa3} and (3.12), we have
\begin{equation}
\pi^{+}_{\xi_{n}}\partial_{x_{n}}\sigma_{-1}(\widetilde{D}_{F}^{-1})(x_{0})|_{|\xi'|=1}
=\frac{\partial_{x_{n}}[c(\xi')](x_{0})}{2(\xi_{n}-\sqrt{-1})}
+\sqrt{-1}h'(0)\bigg[\frac{\sqrt{-1}c(\xi')}{4(\xi_{n}-\sqrt{-1})}
+\frac{c(\xi')+\sqrt{-1}c(dx_{n})}{4(\xi_{n}-\sqrt{-1})^{2}}\bigg].
\end{equation}
By (3.29) and direct calculations, we have\\
\begin{equation}
\partial_{\xi_{n}}\sigma_{-3}((\widetilde{D}^{*}_{F}\widetilde{D}_{F}\widetilde{D}^{*}_{F})^{-1})
=\frac{-4\sqrt{-1}\xi_nc(\xi')+\sqrt{-1}(1-3\xi_n^{2})c(dx_n)}{(1+\xi_n^{2})^3}
\end{equation}
and
\begin{equation}
\partial^{2}_{\xi_{n}}\sigma_{-3}((\widetilde{D}^{*}_{F}\widetilde{D}_{F}\widetilde{D}^{*}_{F})^{-1})
=\sqrt{-1}\bigg[\frac{(20\xi^{2}_{n}-4)c(\xi')+
12(\xi^{3}_{n}-\xi_{n})c(dx_{n})}{(1+\xi_{n}^{2})^{4}}\bigg].
\end{equation}

Since $n=6$, ${\rm trace}_{S(TM)\otimes F}[-{\rm id}]=-8{\rm dim}F$. By the relation of the Clifford action and ${\rm trace}PQ={\rm trace}QP$,  then
\begin{eqnarray}
&&{\rm trace}[c(\xi')c(dx_{n})]=0; \ {\rm trace}[c(dx_{n})^{2}]=-8{\rm dim}F;\
{\rm trace}[c(\xi')^{2}](x_{0})|_{|\xi'|=1}=-8{\rm dim}F;\nonumber\\
&&{\rm trace}[\partial_{x_{n}}[c(\xi')]c(dx_{n})]=0; \
{\rm trace}[\partial_{x_{n}}c(\xi')c(\xi')](x_{0})|_{|\xi'|=1}=-4h'(0){\rm dim}F.
\end{eqnarray}
By (3.41)-(3.44), we get
\begin{eqnarray}
&&{\rm
trace} \Big[\partial_{x_{n}}\pi^{+}_{\xi_{n}}\sigma_{-1}(\widetilde{D}_{F}^{-1})
      \times\partial^{2}_{\xi_{n}}\sigma_{-3}((\widetilde{D}^{*}_{F}\widetilde{D}_{F}
      \widetilde{D}^{*}_{F})^{-1})\Big](x_0)\nonumber\\
&=&h'(0){\rm dim}F\frac{-8-24\xi_{n}\sqrt{-1}+40\xi^{2}_{n}+24\sqrt{-1}\xi^{3}_{n}}{(\xi_{n}-\sqrt{-1})^{6}
(\xi_{n}+\sqrt{-1})^{4}}.
\end{eqnarray}
Then we obtain

\begin{eqnarray}
&&-\frac{1}{2}\int_{|\xi'|=1}\int^{+\infty}_{-\infty}{\rm trace}
\Big[\partial_{x_{n}}\pi^{+}_{\xi_{n}}\sigma_{-1}(\widetilde{D}_{F}^{-1})
      \times\partial^2_{\xi_{n}}\sigma_{-3}
(\widetilde{D}_{F}^{*}\widetilde{D}_{F}\widetilde{D}^{*}_{F})^{-1}\Big](x_0)d\xi_n\sigma(\xi')dx'\nonumber\\
      &=&-\frac{15}{16}\pi h'(0)\Omega_{4}{\rm dim}Fdx'.
\end{eqnarray}

On the other hand, by calculations, we have
\begin{equation}
\pi^{+}_{\xi_{n}}\sigma_{-1}(\widetilde{D}_{F}^{-1})(x_{0})|_{|\xi'|=1}
=-\frac{c(\xi')+\sqrt{-1}c(dx_{n})}{2(\xi_{n}-\sqrt{-1})}.
\end{equation}

By (3.42), (3.44) and (3.47), we get
\begin{equation}
{\rm
trace} \Big[\pi^{+}_{\xi_{n}}\sigma_{-1}(\widetilde{D}_{F}^{-1})
      \times\partial^{2}_{\xi_{n}}\sigma_{-3}((\widetilde{D}^{*}_{F}
      \widetilde{D}_{F}\widetilde{D}^{*}_{F})^{-1})\Big](x_0)
=-16{\rm dim}F\frac{5\xi^2_{n}\sqrt{-1}-\sqrt{-1}-3\xi^{3}_{n}+3\xi_{n}}
{(\xi_{n}-\sqrt{-1})^{5}(\xi_{n}+\sqrt{-1})^{4}}.
\end{equation}
Then we obtain
\begin{eqnarray}
&&-\frac{1}{2}f^{-1}\partial_{x_{n}}(f)
      \int_{|\xi'|=1}\int^{+\infty}_{-\infty}{\rm trace}
\Big[\pi^{+}_{\xi_{n}}\sigma_{-1}(\widetilde{D}_{F}^{-1})
      \times\partial^2_{\xi_{n}}\sigma_{-3}
(\widetilde{D}_{F}^{*}\widetilde{D}_{F}\widetilde{D}^{*}_{F})^{-1}\Big](x_0) d\xi_n\sigma(\xi')dx'\nonumber\\
      &=&\frac{5\sqrt{-1}+44}{4}\pi f^{-1}\partial_{x_{n}}(f)\cdot{\rm dim}F{\rm \Omega_{4}}dx',
\end{eqnarray}
where ${\rm \Omega_{4}}$ is the canonical volume of $S_{4}.$\\
Combining (3.40), (3.46) and (3.49), we obtain
\begin{eqnarray}
{\bf case~(a)~II)}&=&-\frac{15}{16}\pi h'(0)\Omega_{4}{\rm dim}Fdx'+
\frac{5\sqrt{-1}+44}{4}\pi f^{-1}\partial_{x_{n}}(f)\cdot{\rm \Omega_{4}}{\rm dim}Fdx'.
\end{eqnarray}

\noindent  {\bf case (a)~(III)}~$r=-1,l=-3,|\alpha|=j=0,k=1$.\\

By (3.18), we have
 \begin{eqnarray}
{\rm case~ (a)~(III)}&=&-\frac{1}{2}\int_{|\xi'|=1}\int^{+\infty}_{-\infty}{\rm trace} \Big[\partial_{\xi_{n}}\pi^{+}_{\xi_{n}}\sigma_{-1}(f\widetilde{D}_{F}^{-1})
      \times\partial_{\xi_{n}}\partial_{x_{n}}
      \sigma_{-3}
\big(f^{-1}(\widetilde{D}_{F}^{*})^{-1}\cdot f\widetilde{D}_{F}^{-1}\cdot\nonumber\\
&&f^{-1}(\widetilde{D}^{*}_{F})^{-1}\big)\Big](x_0)d\xi_n\sigma(\xi')dx'\nonumber\\
  &=&-\frac{1}{2}\int_{|\xi'|=1}\int^{+\infty}_{-\infty}{\rm trace}
\Big[\partial_{\xi_{n}}\pi^{+}_{\xi_{n}}\big(\sigma_{-1}(\widetilde{D}_{F}^{-1})\big)
      \times\partial_{\xi_{n}}\partial_{x_{n}}\sigma_{-3}
(\widetilde{D}_{F}^{*}\widetilde{D}_{F}\widetilde{D}^{*}_{F})^{-1}\Big](x_0)d\xi_n\sigma(\xi')dx'\nonumber\\
      &&-\frac{1}{2}f\partial_{x_{n}}(f^{-1})
      \int_{|\xi'|=1}\int^{+\infty}_{-\infty}{\rm trace}
\Big[\partial_{\xi_{n}}\pi^{+}_{\xi_{n}}\sigma_{-1}(\widetilde{D}_{F}^{-1})
      \times\partial_{\xi_{n}}\sigma_{-3}
(\widetilde{D}_{F}^{*}\widetilde{D}_{F}\widetilde{D}^{*}_{F})^{-1}\Big](x_0) d\xi_n\sigma(\xi')\nonumber\\
      &&\times dx'.
\end{eqnarray}
By (2.2.29) in \cite{Wa3}, we have
\begin{equation}
\partial_{\xi_{n}}\pi^{+}_{\xi_{n}}\sigma_{-1}(\widetilde{D}_{F}^{-1})(x_{0})|_{|\xi'|=1}
=-\frac{c(\xi')+\sqrt{-1}c(dx_{n})}{2(\xi_{n}-\sqrt{-1})^{2}}.
\end{equation}
By (3.29) and direct calculations, we have\\
\begin{equation}
\partial_{\xi_{n}}\partial_{x_{n}}\sigma_{-3}((\widetilde{D}^{*}_{F}\widetilde{D}_{F}\widetilde{D}^{*}_{F})^{-1})=-\frac{4 \sqrt{-1}\xi_{n}\partial_{x_{n}}c(\xi')(x_{0})}{(1+\xi_{n}^{2})^{3}}
      +\frac{12\sqrt{-1}h'(0)\xi_{n}c(\xi')}{(1+\xi_{n}^{2})^{4}}
      -\frac{\sqrt{-1}(2-10\xi^{2}_{n})h'(0)c(dx_{n})}{(1+\xi_{n}^{2})^{4}}.
\end{equation}

Combining (3.44), (3.52) and (3.53), we have
\begin{eqnarray}
&&{\rm trace} \Big[\partial_{\xi_{n}}\pi^{+}_{\xi_{n}}\sigma_{-1}(\widetilde{D}_{F}^{-1})
      \times\partial_{\xi_{n}}\partial_{x_{n}}\sigma_{-3}((\widetilde{D}^{*}_{F}
      \widetilde{D}_{F}\widetilde{D}^{*}_{F})^{-1})\Big](x_{0})|_{|\xi'|=1}\nonumber\\
&=&h'(0){\rm dim}F\frac{8\sqrt{-1}-32\xi_{n}-8\sqrt{-1}\xi^{2}_{n}}{(\xi_{n}-\sqrt{-1})^{5}(\xi+i)^{4}},
\end{eqnarray}
and
\begin{eqnarray}
&&{\rm trace} \Big[\partial_{\xi_{n}}\pi^{+}_{\xi_{n}}\sigma_{-1}(\widetilde{D}_{F}^{-1})
      \times\partial_{\xi_{n}}\sigma_{-3}((\widetilde{D}^{*}_{F}\widetilde{D}_{F}
      \widetilde{D}^{*}_{F})^{-1})\Big](x_{0})|_{|\xi'|=1}\nonumber\\
&=&-4{\rm dim}F\frac{4\sqrt{-1}\xi_n+1-3\xi^2_{n}}{(\xi_{n}-\sqrt{-1})^{5}(\xi_n+\sqrt{-1})^{3}}.
\end{eqnarray}
Then
\begin{eqnarray}
&&-\frac{1}{2}\int_{|\xi'|=1}\int^{+\infty}_{-\infty}{\rm trace}
\Big[\partial_{x_{n}}\pi^{+}_{\xi_{n}}\big(\sigma_{-1}(\widetilde{D}_{F}^{-1})\big)
      \times\partial_{\xi_{n}}\partial_{x_{n}}\sigma_{-3}
(\widetilde{D}_{F}^{*}\widetilde{D}_{F}\widetilde{D}^{*}_{F})^{-1}\Big](x_0)d\xi_n\sigma(\xi')dx'\nonumber\\
      &=&\frac{25}{16}\pi h'(0)\Omega_{4}{\rm dim}Fdx',
\end{eqnarray}
and
\begin{eqnarray}
&&-\frac{1}{2}f\partial_{x_{n}}(f^{-1})
      \int_{|\xi'|=1}\int^{+\infty}_{-\infty}{\rm trace}
\Big[\partial_{\xi_{n}}\pi^{+}_{\xi_{n}}\sigma_{-1}(\widetilde{D}_{F}^{-1})
      \times\partial_{\xi_{n}}\sigma_{-3}
(\widetilde{D}_{F}^{*}\widetilde{D}_{F}\widetilde{D}^{*}_{F})^{-1}\Big](x_0) d\xi_n\sigma(\xi')dx'\nonumber\\
      &=&\frac{\pi \sqrt{-1}}{16}\cdot f\cdot\partial_{x_{n}}(f^{-1})\Omega_{4}{\rm dim}Fdx',
\end{eqnarray}
where ${\rm \Omega_{4}}$ is the canonical volume of $S_{4}.$

Then
\begin{eqnarray}
{\bf case~(a)~III)}&=&\Big[\frac{25}{16}\pi h'(0)+
\frac{\pi \sqrt{-1}}{16}\cdot f\cdot\partial_{x_{n}}(f^{-1})\Big]\Omega_{4}{\rm dim}Fdx'.
\end{eqnarray}
\\
\noindent  {\bf case (b)}~$r=-1,l=-4,|\alpha|=j=k=0$.\\

By (3.18), we have
 \begin{eqnarray}
{\rm case~ (b)}&=&-i\int_{|\xi'|=1}\int^{+\infty}_{-\infty}{\rm trace} \Big[\pi^{+}_{\xi_{n}}\sigma_{-1}(f\widetilde{D}_{F}^{-1})
      \times\partial_{\xi_{n}}\sigma_{-4}
      \big(f^{-1}(\widetilde{D}_{F}^{*})^{-1}\cdot f\widetilde{D}_{F}^{-1}\cdot f^{-1}(\widetilde{D}^{*}_{F})^{-1}\big)\Big](x_0)\nonumber\\
&&\times d\xi_n\sigma(\xi')dx'\nonumber\\
&=&-i\int_{|\xi'|=1}\int^{+\infty}_{-\infty}{\rm trace} \Bigg[\pi^{+}_{\xi_{n}}\sigma_{-1}(f\widetilde{D}_{F}^{-1})
      \times\partial_{\xi_{n}}\Bigg(f^{-1}\sigma_{-4}(\widetilde{D}^{*}_{F}
      \widetilde{D}_{F}\widetilde{D}^{*}_{F})^{-1}+\frac{2c(\xi)c(df)c(\xi)}{f^{2}|\xi|^6}
\nonumber\\
      &&+\frac{ic(\xi)\sum\limits_j\Big[c(dx_j)|\xi|^2+2\xi_{j}c(\xi)\Big]D_{x_j}(f^{-1})c(\xi)}{|\xi|^8}
      \Bigg)\Bigg](x_0)d\xi_n\sigma(\xi')dx'\nonumber\\
     &=&-i\int_{|\xi'|=1}\int^{+\infty}_{-\infty}{\rm trace} \Big[\pi^{+}_{\xi_{n}}\sigma_{-1}(\widetilde{D}_{F}^{-1})
      \times\partial_{\xi_{n}}\Big(\sigma_{-4}(\widetilde{D}^{*}_{F}
      \widetilde{D}_{F}\widetilde{D}^{*}_{F})^{-1}\Big)\Big](x_0)d\xi_n\sigma(\xi')dx'\nonumber\\
&&-2if^{-1}\int_{|\xi'|=1}\int^{+\infty}_{-\infty}{\rm trace} \Bigg[\pi^{+}_{\xi_{n}}\sigma_{-1}(\widetilde{D}_{F}^{-1})
      \times\partial_{\xi_{n}}\Bigg(\frac{c(\xi)c(df)c(\xi)}{|\xi|^6}
      \Bigg)\Bigg](x_0)d\xi_n\sigma(\xi')dx'\nonumber\\
      &&-fi\int_{|\xi'|=1}\int^{+\infty}_{-\infty}{\rm trace} \Big[\pi^{+}_{\xi_{n}}\sigma_{-1}(\widetilde{D}_{F}^{-1})
      \times\partial_{\xi_{n}}\Big(\frac{ic(\xi)\sum\limits_j\Big[c(dx_j)|\xi|^2
      +2\xi_{j}c(\xi)\Big]D_{x_j}(f^{-1})c(\xi)}{|\xi|^8}\Big)\Big]\nonumber\\
      &&\times(x_0)d\xi_n\sigma(\xi')dx',
\end{eqnarray}

In the normal coordinate, $g^{ij}(x_{0})=\delta^{j}_{i}$ and $\partial_{x_{j}}(g^{\alpha\beta})(x_{0})=0$, if $j<n$; $\partial_{x_{j}}(g^{\alpha\beta})(x_{0})=h'(0)\delta^{\alpha}_{\beta}$, if $j=n$.
So by Lemma A.2 in \cite{Wa3}, we have $\Gamma^{n}(x_{0})=\frac{5}{2}h'(0)$ and $\Gamma^{k}(x_{0})=0$ for $k<n$. By the definition of $\delta^{k}$ and Lemma 2.3 in \cite{Wa3}, we have $\delta^{n}(x_{0})=0$ and $\delta^{k}=\frac{1}{4}h'(0)c(\widetilde{e}_{k})c(\widetilde{e}_{n})$ for $k<n$. By (3.30), we obtain
\begin{eqnarray}
&&\sigma_{-4}(\widetilde{D}^{*}_{F}\widetilde{D}_{F}\widetilde{D}^{*}_{F})^{-1}\nonumber\\
&=&\frac{-17-9\xi^{2}_{n}}{4(1+\xi^{2}_{n})^{4}}h'(0)c(\xi')c(dx_{n})c(\xi')+\frac{33\xi_{n}+17\xi^{3}_{n}}{2(1+\xi^{2}_{n})^{4}}h'(0)c(\xi')
+\frac{49\xi^{2}_{n}+25\xi^{4}_{n}}{2(1+\xi^{2}_{n})^{4}}h'(0)c(dx_{n})\nonumber\\
&&+\frac{1}{(1+\xi^{2}_{n})^{3}}c(\xi')c(dx_{n})\partial_{x_{n}}[c(\xi')](x_{0})-\frac{3\xi_{n}}{(1+\xi^{2}_{n})^{3}}\partial_{x_{n}}[c(\xi')](x_{0})
-\frac{2\xi_{n}}{(1+\xi^{2}_{n})^{3}}h'(0)\xi_{n}c(\xi')(x_{0})\nonumber\\
&&+\frac{1-\xi^{2}_{n}}{(1+\xi^{2}_{n})^{3}}h'(0)c(dx_{n})(x_{0})+\frac{c(\xi)\lambda c(\xi)}{|\xi|^6}-\frac{2c(\xi)c(A^*)c(\xi)}{|\xi|^6}
-\frac{2c(A)}{|\xi|^4}.
\end{eqnarray}
Then
\begin{eqnarray}
&&\partial_{\xi_{n}}\Big(\sigma_{-4}(\widetilde{D}^{*}_{F}\widetilde{D}_{F}
\widetilde{D}^{*}_{F})^{-1}\Big)(x_{0})\nonumber\\
&=&\frac{59\xi_{n}+27\xi^{3}_{n}}{2(1+\xi^{2}_{n})^{5}}h'(0)c(\xi')c(dx_{n})c(\xi')
+\frac{33-180\xi^{2}_{n}-85\xi^{4}_{n}}{2(1+\xi^{2}_{n})^{5}}h'(0)c(\xi')+\frac{49\xi_{n}-97\xi^{3}_{n}
-50\xi^{5}_{n}}{2(1+\xi^{2}_{n})^{5}}h'(0)c(dx_{n})\nonumber\\
&&
-\frac{6\xi_{n}}{(1+\xi^{2}_{n})^{4}}c(\xi')c(dx_{n})\partial_{x_{n}}[c(\xi')](x_{0})-\frac{3-15\xi^{2}_{n}}{(1+\xi^{2}_{n})^{4}}\partial_{x_{n}}[c(\xi')](x_{0})
+\frac{4\xi_{n}^3-8\xi_{n}}{(1+\xi^{2}_{n})^{4}}h'(0)c(dx_{n})\nonumber\\
&&
+\frac{2-10\xi^2_{n}}{(1+\xi^{2}_{n})^{4}}h'(0)c(\xi')+\frac{c(dx_{n})\lambda c(\xi')+c(\xi')\lambda c(dx_{n})+2\xi_{n}c(dx_{n})\lambda c(dx_{n})}{(1+\xi^{2}_{n})^{3}}-\frac{6\xi_{n}c(\xi)\lambda c(\xi)}{(1+\xi^{2}_{n})^{4}}\nonumber\\
&&+\frac{c(dx_{n})c(A^*) c(\xi')+c(\xi')c(A^*) c(dx_{n})+2\xi_{n}c(dx_{n})c(A^*)c(dx_{n})}{(1+\xi^{2}_{n})^{3}}-\frac{6\xi_{n}c(\xi)c(A^*) c(\xi)}{(1+\xi^{2}_{n})^{4}}\nonumber\\
&&-\frac{2\xi_{n}c(A)}{(1+\xi^{2}_{n})^{3}}.
\end{eqnarray}

By (3.47) and (3.61), we obtain
\begin{eqnarray}
&&{\rm trace} \Big[\pi^{+}_{\xi_{n}}\sigma_{-1}(\widetilde{D}_{F}^{-1})
      \times\partial_{\xi_{n}}\sigma_{-4}(\widetilde{D}^{*}_{F}\widetilde{D}_{F}
      \widetilde{D}^{*}_{F})^{-1}\Big](x_{0})|_{|\xi'|=1}\nonumber\\
&=&h'(0){\rm dim}F\frac{4i(-17-42i\xi_{n}+50\xi^{2}_{n}
-16i\xi^{3}_{n}+29\xi^{4}_{n})}{(\xi_{n}-i)^{5}(\xi+i)^{5}}\nonumber\\
&&+\frac{(4\xi_{n}i+2)i}{2(\xi_{n}+i)(1+\xi^{2}_{n})^{3}}{\rm trace}[c(\xi')\lambda]+\frac{4\xi_{n}i+2}{2(\xi_{n}+i)(1+\xi^{2}_{n})^{3}}{\rm trace}[c(dx_{n})\lambda]\nonumber\\
&&+\frac{(4\xi_{n}i+2)i}{2(\xi_{n}+i)(1+\xi^{2}_{n})^{3}}{\rm trace}[c(\xi')c(A^*)]+\frac{4\xi_{n}i+2}{2(\xi_{n}+i)(1+\xi^{2}_{n})^{3}}{\rm trace}[c(dx_{n})c(A^*)]\nonumber\\
&&+\frac{-2\xi_{n}}{2(\xi_{n}-i)(1+\xi^{2}_{n})^{3}}{\rm trace}[c(\xi')c(A)]+\frac{-2\xi_{n}i}{2(\xi_{n}-i)(1+\xi^{2}_{n})^{3}}{\rm trace}[c(dx_{n})c(A)].
\end{eqnarray}
By the relation of the Clifford action and ${\rm trace}QP={\rm trace}PQ$, then we have the following equalities
\begin{eqnarray}
&&{\rm trace}\Big[c(dx_{n})\lambda\Big]={\rm trace}\Big[c(dx_{n})\sum_{j=1}^{n}c(e_{j})(\sigma_{j}^{F}-A^{*}(e_{j}))\Big]={\rm trace}\Big[-{\rm id}\otimes
(\sigma_{n}^{F}-A^{*}(e_{n}))\Big];\\
&&{\rm trace}\Big[c(\xi')\lambda\Big]={\rm trace}\Big[c(\xi')\sum_{j=1}^{n}c(e_{j})(\sigma_{j}^{F}-A^{*}(e_{j}))\Big]={\rm trace}\Big[-\sum_{j=1}^{n-1}\xi_j
(\sigma_{j}^{F}-A^{*}(e_{j}))\Big];\\
&&{\rm trace}\Big[c(dx_{n})c(A^*)\Big]={\rm trace}\Big[c(dx_{n})\sum_{j=1}^{n}c(e_{j})\otimes A^*(e_{j})\Big]={\rm trace}\Big[-{\rm id}\otimes A^*(e_{n})\Big];\\
&&{\rm trace}\Big[c(dx_{n})c(A)\Big]={\rm trace}\Big[c(dx_{n})\sum_{j=1}^{n}c(e_{j})\otimes A(e_{j})\Big]={\rm trace}\Big[-{\rm id}\otimes A(e_{n})\Big];\\
&&{\rm trace}\Big[c(\xi')c(A^*)\Big]={\rm trace}\Big[c(\xi')\sum_{j=1}^{n}c(e_{j})\otimes A^*(e_{j})\Big]={\rm trace}\Big[-\sum_{j=1}^{n-1}\xi_jA^*(e_{j})\Big];\\
&&{\rm trace}\Big[c(\xi')c(A)\Big]={\rm trace}\Big[c(\xi')\sum_{j=1}^{n}c(e_{j})\otimes A(e_{j})\Big]={\rm trace}\Big[-\sum_{j=1}^{n-1}\xi_jA(e_{j})\Big].
\end{eqnarray}
We note that $i<n,~\int_{|\xi'|=1}\xi_i\sigma(\xi')=0$,
so ${\rm trace}\big[c(\xi')c(A^*)\big]$ has no contribution for computing {\bf case (b)}.

By (3.24), then
 \begin{eqnarray}
&&-i\int_{|\xi'|=1}\int^{+\infty}_{-\infty}{\rm trace} \Big[\pi^{+}_{\xi_{n}}\sigma_{-1}(\widetilde{D}_{F}^{-1})
      \times\partial_{\xi_{n}}\Big(\sigma_{-4}(\widetilde{D}^{*}_{F}
      \widetilde{D}_{F}\widetilde{D}^{*}_{F})^{-1}\Big)\Big](x_0)d\xi_n\sigma(\xi')dx'\nonumber\\
&=&\Bigg\{-\frac{129}{16} h'(0)+\frac{3}{2}{\rm trace}\Big[\sigma_{n}^{F}-A^{*}(e_{n})\Big]-3{\rm trace}\Big[A^{*}(e_{n})\Big]- {\rm trace}\Big[A(e_{n})\Big]\Bigg\}\pi {\rm dim}F\Omega_{4}dx'.
\end{eqnarray}
Since
\begin{eqnarray*}
\partial_{\xi_{n}}\Big(\frac{c(\xi)c(df) c(\xi)}{|\xi|^{6}}\Big)&=&\frac{c(dx_{n})c(df) c(\xi')+c(\xi')c(df) c(dx_{n})+2\xi_{n}c(dx_{n})c(df) c(dx_{n})}{(1+\xi^{2}_{n})^{3}}\nonumber\\
&&-\frac{6\xi_{n}c(\xi)c(df) c(\xi)}{(1+\xi^{2}_{n})^{4}}
\end{eqnarray*}
and
\begin{eqnarray}
&&\partial_{\xi_{n}}\Big(\frac{ic(\xi)\sum\limits_j\Big[c(dx_j)|\xi|^2
+2\xi_{j}c(\xi)\Big]D_{x_j}(f^{-1})c(\xi)}{|\xi|^8}\Big)\nonumber\\
&=&i\Bigg\{c(dx_{n})\sum\limits_j\Big[c(dx_j)|\xi|^2
+2\xi_{j}c(\xi)\Big]D_{x_j}(f^{-1}) c(\xi')
+c(\xi')\sum\limits_j\Big[c(dx_j)|\xi|^2
+2\xi_{j}c(\xi)\Big]D_{x_j}(f^{-1}) c(dx_{n})\nonumber\\
&&+2\xi_{n}c(dx_{n})\sum\limits_j\Big[c(dx_j)|\xi|^2
+2\xi_{j}c(\xi)\Big]D_{x_j}(f^{-1}) c(dx_{n})\Bigg\}(1+\xi^{2}_{n})^{-4}-i\Bigg\{8\xi_{n}c(\xi)\sum\limits_j\Big[c(dx_j)|\xi|^2
\nonumber\\
&&+2\xi_{j}c(\xi)\Big]D_{x_j}(f^{-1}) c(\xi)\Bigg\}(1+\xi^{2}_{n})^{-5},
\end{eqnarray}
then
\begin{eqnarray*}
&&{\rm trace} \Big[\pi^{+}_{\xi_{n}}\sigma_{-1}(\widetilde{D}_{F}^{-1})
      \times\partial_{\xi_{n}}\Big(\frac{c(\xi)c(df) c(\xi)}{|\xi|^{6}}\Big)\Big](x_0) \nonumber\\
&=&\frac{(4\xi_{n}i+2)i}{2(\xi_{n}+i)(1+\xi^{2}_{n})^{3}}{\rm trace}[c(\xi')c(df)]+\frac{4\xi_{n}i+2}{2(\xi_{n}+i)(1+\xi^{2}_{n})^{3}}{\rm trace}[c(dx_{n})c(df)].
\end{eqnarray*}
and
\begin{eqnarray}
&&{\rm trace} \Bigg[\pi^{+}_{\xi_{n}}\sigma_{-1}(\widetilde{D}_{F}^{-1})
      \times\partial_{\xi_{n}}\Big(\frac{ic(\xi)\sum\limits_j\Big[c(dx_j)|\xi|^2
+2\xi_{j}c(\xi)\Big]D_{x_j}(f^{-1})c(\xi)}{|\xi|^8}\Big)\Bigg](x_0) \nonumber\\
&=&\frac{(3\xi_{n}-i)i}{(\xi_{n}+i)(1+\xi^{2}_{n})^{4}}{\rm trace}\Bigg[c(\xi')\sum\limits_j\Big[c(dx_j)|\xi|^2
+2\xi_{j}c(\xi)\Big]D_{x_j}(f^{-1})\Bigg]\nonumber\\
&&+\frac{3\xi_{n}-i}{(\xi_{n}+i)(1+\xi^{2}_{n})^{4}}{\rm trace}\Bigg[c(dx_{n})\sum\limits_j\Big[c(dx_j)|\xi|^2
+2\xi_{j}c(\xi)\Big]D_{x_j}(f^{-1})\Bigg].
\end{eqnarray}

By the relation of the Clifford action and ${\rm trace}QP={\rm trace}PQ$, then we have the following equalities
\begin{eqnarray*}
&&{\rm trace}\Big[c(dx_{n})c(df)\Big]=-g(dx_{n},df);\\
\end{eqnarray*}
and
\begin{eqnarray*}
&&{\rm trace}\Bigg[c(dx_{n})\sum\limits_j\Big[c(dx_j)|\xi|^2
      +2\xi_{j}c(\xi)\Big]D_{x_j}(f^{-1})\Bigg]\nonumber\\
      &=&{\rm trace}\big(-{\rm id}\big)|\xi|^2\bigg(-i\partial_{x_n}(f)f^{-1}\bigg)
      +2\sum\limits_{j}\xi_j\xi_n{\rm trace}\big(-{\rm id}\big)\bigg(-i\partial_{x_j}(f)f^{-1}\bigg)\nonumber\\
      &=&-8{\rm dim}F|\xi|^2\bigg(-i\partial_{x_n}(f)f^{-1}\bigg)+2\sum\limits_{j}\xi_j\xi_n{\rm trace}\big(-{\rm id}\big)\bigg(-i\partial_{x_j}(f)f^{-1}\bigg).
\end{eqnarray*}
We note that $i<n,~\int_{|\xi'|=1}\xi_i\sigma(\xi')=0$,
so ${\rm trace}\big[c(\xi')c(df)\big]$, ${\rm trace}\Bigg[c(\xi')\sum\limits_j\Big[c(dx_j)|\xi|^2
      +2\xi_{j}c(\xi)\Big]D_{x_j}(f^{-1})\Bigg]$ and $2i\sum\limits_{j}\xi_j\xi_n\partial_{x_j}(f)f^{-1}{\rm trace}[-{\rm id}]$ have no contribution for computing {\bf case (b)}.
Then we obtain
\begin{eqnarray}
&&-2if^{-1}\int_{|\xi'|=1}\int^{+\infty}_{-\infty}{\rm trace} \Bigg[\pi^{+}_{\xi_{n}}\sigma_{-1}(\widetilde{D}_{F}^{-1})
      \times\partial_{\xi_{n}}\Big(\frac{c(\xi)c(df) c(\xi)}{|\xi|^{6}}\Big)\Bigg](x_0)d\xi_n\sigma(\xi')dx'\nonumber\\
      &=&\frac{3}{8f}\pi  g(dx_{n},df)\Omega_{4}dx'.
\end{eqnarray}
and
\begin{eqnarray}
&&-fi\int_{|\xi'|=1}\int^{+\infty}_{-\infty}{\rm trace} \Big[\pi^{+}_{\xi_{n}}\sigma_{-1}(\widetilde{D}_{F}^{-1})
      \times\partial_{\xi_{n}}\Big(\frac{ic(\xi)\sum\limits_j\Big[c(dx_j)|\xi|^2
      +2\xi_{j}c(\xi)\Big]D_{x_j}(f^{-1})c(\xi)}{|\xi|^8}\Big)\Big]\nonumber\\
      &&\times(x_0)d\xi_n\sigma(\xi')dx'\nonumber\\
      &=&-\frac{15i}{2}\partial_{x_{n}}(f)\pi  {\rm dim}F\Omega_{4}dx'.
\end{eqnarray}

Thus we have
\begin{eqnarray}
{\bf case~ (b)}
&=&\Big\{-\frac{129}{16} h'(0)+\frac{3}{2} {\rm trace}\Big[\sigma_{n}^{F}-A^{*}(e_{n})\Big]-3 {\rm trace}\Big[A^{*}(e_{n})\Big]-{\rm trace}\Big[A(e_{n})\Big]\Big\}\pi {\rm dim}F \Omega_{4}dx'\nonumber\\
&&+\frac{3}{8f}\pi  g[dx_{n},df]\Omega_{4}dx'
      -\frac{15i}{2}\partial_{x_{n}}(f)\pi  {\rm dim}F\Omega_{4}dx'.
\end{eqnarray}

\noindent {\bf  case (c)}~$r=-2,l=-3,|\alpha|=j=k=0$.\\

By (3.18), we have

\begin{eqnarray}
{\rm case~ (c)}&=&-i\int_{|\xi'|=1}\int^{+\infty}_{-\infty}{\rm trace} \Big[\pi^{+}_{\xi_{n}}\sigma_{-2}(f\widetilde{D}_{F}^{-1})
      \times\partial_{\xi_{n}}
      \sigma_{-3}
\big(f^{-1}(\widetilde{D}_{F}^{*})^{-1}\cdot f\widetilde{D}_{F}^{-1}\cdot f^{-1}(\widetilde{D}^{*}_{F})^{-1}\big)\Big](x_0)\nonumber\\
&&\times d\xi_n\sigma(\xi')dx'\nonumber\\
&=&-i\int_{|\xi'|=1}\int^{+\infty}_{-\infty}{\rm trace} \Big[\pi^{+}_{\xi_{n}}\sigma_{-2}(\widetilde{D}_{F}^{-1})
      \times\partial_{\xi_{n}}
      \sigma_{-3}
\big((\widetilde{D}_{F}^{*}\widetilde{D}_{F}\widetilde{D}^{*}_{F})^{-1}\big)\Big]
(x_0)d\xi_n\sigma(\xi')dx'.
\end{eqnarray}

By (3.14), we have
\begin{eqnarray}
\pi^{+}_{\xi_{n}}\sigma_{-2}(\widetilde{D}_{F}^{-1})
&=&\pi^{+}_{\xi_{n}}\Big(\frac{c(\xi)\sigma_0(\widetilde{D}_{F})c(\xi)}{|\xi|^4}+\frac{c(\xi)}{|\xi|^6}\sum_jc(dx_j)
\Big[\partial_{x_j}[c(\xi)]|\xi|^2-c(\xi)\partial_{x_j}(|\xi|^2)\Big]\Big)\nonumber\\
&:=&T_{1}-T_{2}+\pi^{+}_{\xi_{n}}\Big(\frac{c(\xi)\mu c(\xi)}{|\xi|^{4}}\Big),
\end{eqnarray}

where
\begin{eqnarray}
T_{1}&=&-\frac{1}{4(\xi_{n}-i)^{2}}\Big[(2+i\xi_{n})c(\xi')\sigma_{0}(\widetilde{D}_{F})c(\xi')
+i\xi_{n}c(dx_{n})\sigma_{0}(\widetilde{D}_{F})c(dx_{n})+(2+i\xi_{n})c(\xi')c(dx_{n})\nonumber\\
&&\times\partial_{x_{n}}[c(\xi')]+ic(dx_{n})\sigma_{0}(\widetilde{D}_{F})c(\xi')+ic(\xi')
\sigma_{0}(\widetilde{D}_{F})c(dx_{n})-i\partial_{x_{n}}[c(\xi')]  \Big]\nonumber\\
&=&\frac{1}{4(\xi_{n}-i)^{2}}\Big[\frac{5}{2}h'(0)c(dx_{n})-\frac{5i}{2}h'(0)c(\xi')-(2+i\xi_{n})c(\xi')c(dx_{n})\partial_{\xi_{n}}[c(\xi')]
+i\partial_{\xi_{n}}[c(\xi')] \Big];\\
T_{2}&=&\frac{h'(0)}{2}\Big[\frac{c(dx_{n})}{4i(\xi_{n}-i)}+\frac{c(dx_{n})-ic(\xi')}{8(\xi_{n}-i)^{2}}
+\frac{3\xi_{n}-7i}{8(\xi_{n}-i)^{3}}\big(ic(\xi')-c(dx_{n})\big)\Big].
\end{eqnarray}

On the other hand,
\begin{eqnarray}
&&\pi^+_{\xi_n}\Big( \frac{c(\xi)\mu c(\xi)}{|\xi|^4}\Big)(x_{0})|_{|\xi'|=1}\nonumber\\
 &=&\frac{(-i\xi_n-2)c(\xi')\mu c(\xi')-i\Big[c(dx_{n})\mu c(\xi')+c(\xi')\mu c(dx_{n})\Big]
  -i\xi_n c(dx_{n})\mu c(dx_{n})  }{4(\xi_n-i)^{2}}.
\end{eqnarray}

By (3.42) (3.44) and (3.76), then we have
\begin{eqnarray}
&&{\rm tr}
\Big[T_1\times\partial_{\xi_n}\sigma_{-3}((\widetilde{D}^{*}_{F}
\widetilde{D}_{F}\widetilde{D}^{*}_{F})^{-1})\Big]|_{|\xi'|=1}\nonumber\\
&=&{\rm tr }\Big\{ \frac{1}{4(\xi_n-i)^2}\Big[\frac{5}{2}h'(0)c(dx_n)-\frac{5i}{2}h'(0)c(\xi')
  -(2+i\xi_n)c(\xi')c(dx_n)\partial_{\xi_n}c(\xi')+i\partial_{\xi_n}c(\xi')\Big]\nonumber\\
&&\times \frac{-4i\xi_nc(\xi')+(i-3i\xi_n^{2})c(dx_n)}{(1+\xi_n^{2})^3}\Big\} \nonumber\\
&=&h'(0){\rm dim}F\frac{3+12i\xi_n+3\xi_n^{2}}{(\xi_n-i)^4(\xi_n+i)^3}.
\end{eqnarray}

Similarly, we have\begin{eqnarray}
&&{\rm trace }\Big[T_2\times\partial_{\xi_n}\sigma_{-3}((\widetilde{D}^{*}_{F}
\widetilde{D}_{F}\widetilde{D}^{*}_{F})^{-1})\Big]|_{|\xi'|=1}\nonumber\\
&=&{\rm trace }\Big\{ \frac{h'(0)}{2}\Big[\frac{c(dx_n)}{4i(\xi_n-i)}+\frac{c(dx_n)-ic(\xi')}{8(\xi_n-i)^2}
+\frac{3\xi_n-7i}{8(\xi_n-i)^3}\Big(ic(\xi')-c(dx_n)\Big)\Big]\nonumber\\
&&\times \frac{-4i\xi_nc(\xi')+(i-3i\xi_n^{2})c(dx_n)}{(1+\xi_n^{2})^3}\Big\} \nonumber\\
&=&h'(0){\rm dim}F\frac{4i-11\xi_n-6i\xi_n^{2}+3\xi_n^{3}}{(\xi_n-i)^5(\xi_n+i)^3}.
\end{eqnarray}

By (3.79) and (3.80), we obtain
\begin{eqnarray}
&&-i\int_{|\xi'|=1}\int^{+\infty}_{-\infty}{\rm trace} \Bigg[\Big(T_{1}-T_{2}\Big)\times\partial_{\xi_{n}}\sigma_{-3}\Big((\widetilde{D}^{*}_{F}\widetilde{D}_{F}
      \widetilde{D}^{*}_{F})^{-1}\Big)\Bigg](x_0)d\xi_n\sigma(\xi')dx'\nonumber\\
&=&-i{\rm dim}Fh'(0)\int_{|\xi'|=1}\int^{+\infty}_{-\infty}
\frac{-7i+26\xi_n+15i\xi_n^{2}}{(\xi_n-i)^5(\xi_n+i)^3}d\xi_n\sigma(\xi')dx' \nonumber\\
&=&-i{\rm dim}Fh'(0)\frac{2 \pi i}{4!}\Big[\frac{-7i+26\xi_n+15i\xi_n^{2}}{(\xi_n+i)^3}
     \Big]^{(5)}|_{\xi_n=i}\Omega_4dx'\nonumber\\
&=&\frac{55}{16}{\rm dim}F\pi h'(0)\Omega_4dx'.
\end{eqnarray}

By (3.55) and (3.56), we have
\begin{eqnarray}
&&{\rm trace} \Big[\pi^+_{\xi_n}\Big( \frac{c(\xi)\mu c(\xi)}{|\xi|^4}\Big)
      \times\partial_{\xi_{n}}\sigma_{-3}((\widetilde{D}^{*}_{F}\widetilde{D}_{F}\widetilde{D}^{*}_{F})^{-1})\Big](x_0) \nonumber\\
&=&\frac{(3\xi_{n}-i)i}{2(\xi_{n}-i)(1+\xi^{2}_{n})^{3}}{\rm trace}[c(dx_{n})\mu]+\frac{3\xi_{n}-i}{2(\xi_{n}-i)(1+\xi^{2}_{n})^{3}}{\rm trace}[c(\xi')\mu].
\end{eqnarray}

By the relation of the Clifford action and ${\rm trace}PQ={\rm trace}QP$, then we have the equalities
\begin{eqnarray}
&&{\rm trace}\Big[c(dx_{n})\mu\Big]
={\rm trace}\Big[c(dx_{n})\sum_{j=1}^{n}c(e_{j})(\sigma_{j}^{F}+A(e_{j}))\Big]={\rm trace}\Big[-{\rm id}\otimes
(\sigma_{n}^{F}+A(e_{n}))\Big];\\
&&{\rm trace}\Big[c(\xi')\mu\Big]
={\rm trace}\Big[c(\xi')\sum_{j=1}^{n}c(e_{j})(\sigma_{j}^{F}+A(e_{j}))\Big]={\rm trace}\Big[-\sum_{j=1}^{n-1}\xi_j
(\sigma_{j}^{F}+A(e_{j}))\Big].
\end{eqnarray}

We note that $i<n,~\int_{|\xi'|=1}\xi_i\sigma(\xi')=0$,
so ${\rm trace}[c(\xi')\mu]$ has no contribution for computing {\bf  case (c)}.

Then, we obtain
\begin{eqnarray}
&&-i\int_{|\xi'|=1}\int^{+\infty}_{-\infty}{\rm trace} \Big[\pi^+_{\xi_n}\Big( \frac{c(\xi)\mu c(\xi)}{|\xi|^4}\Big)
      \times\partial_{\xi_{n}}\sigma_{-3}\Big((\widetilde{D}^{*}_{F}\widetilde{D}_{F}\widetilde{D}^{*}_{F})^{-1}\Big)\Big](x_0)d\xi_n\sigma(\xi')dx'\nonumber\\
&=&-i\int_{|\xi'|=1}\int^{+\infty}_{-\infty}\frac{(3\xi_{n}-i)i}{2(\xi_{n}-i)(1+\xi^{2}_{n})^{3}}{\rm trace}[c(dx_{n})\mu]d\xi_n\sigma(\xi')dx' \nonumber\\
&=&-2\pi {\rm dim}F{\rm trace}[\sigma_{n}^{F}+A(e_{n})]\Omega_4dx'.
\end{eqnarray}

Then
 \begin{eqnarray}
{\bf case~ (c)}&=&\frac{55}{16}{\rm dim}F\pi h'(0)\Omega_4dx'-2\pi {\rm dim}F h'(0){\rm trace}[\sigma_{n}^{F}+A(e_{n})]\Omega_4dx'.
\end{eqnarray}

Now $\Phi$ is the sum of the {\bf  case (a)}, {\bf  case (b)} and {\bf  case (c)}, then
\begin{eqnarray}
\Phi&=&\Big[-4h'(0)-{\rm trace}\Big(A(e_{n})\Big)-3{\rm trace}\Big(A^{*}(e_{n})\Big)+\frac{3}{2}{\rm trace}\Big(\sigma_{n}^{F}-A^{*}(e_{n})\Big)\nonumber\\
&&-2{\rm trace}\Big(\sigma_{n}^{F}+A(e_{n})\Big)
+(\frac{19i}{16}+11)\cdot f^{-1}\cdot\partial_{x_{n}}(f)\Big]\pi {\rm dim}F\Omega_4dx'\nonumber\\
&&
+\frac{3}{8f}\pi  g(dx_{n},df)\Omega_{4}dx'
      -\frac{15i}{2}\partial_{x_{n}}(f)\pi  {\rm dim}F\Omega_{4}dx'.
\end{eqnarray}

By (4.2) in \cite{Wa3}, we have
$$K=\sum_{1\leq i,j\leq n-1}K_{i,j}g^{i,j}_{\partial M};K_{i,j}=-\Gamma^{n}_{i,j},$$
and $K_{i,j}$ is the second fundamental form, or extrinsic curvature. For $n=6$, then
\begin{eqnarray}
K(x_{0})=\sum_{1\leq i,j\leq n-1}K_{i,j}(x_{0})g^{i,j}_{\partial M}(x_{0})=\sum^{5}_{i=1}K_{i,i}(x_{0})=-\frac{5}{2}h'(0).
\end{eqnarray}
Hence we conclude that
\begin{thm}
 Let M be a 6-dimensional compact spin manifolds with the boundary $\partial M$. Then
 \begin{eqnarray}
&&\widetilde{Wres}[\pi^{+}(f\widetilde{D}_{F}^{-1}) \circ\pi^{+}\big(f^{-1}(\widetilde{D}_{F}^{*})^{-1}\cdot f\widetilde{D}_{F}^{-1}\cdot f^{-1}(\widetilde{D}^{*}_{F})^{-1}\big)]\nonumber\\
&=& 8 \pi^{3}\int_{M}\bigg \{{\rm{trace}}\Big[-\frac{s}{12}+c(A^*)c(A)-\frac{1}{4}\sum_{i}\big[c(A^*)c(e_{i})-c(e_{i})c(A) \big]^{2}-\frac{1}{2}\sum_{j}\nabla_{e_{j}}^{F}\big(c(A^*)\big)c(e_{j})
 \nonumber\\
&&-\frac{1}{2}\sum_{j}c(e_{j})\nabla_{e_{j}}^{F}\big(c(A)\big)\Big]
 -\frac{2\Delta(f)}{f}+\frac{4{\rm{trace}}\Big[A(grad_{M}f)\Big]}{f}-f^{2}\Big[|grad_{M}(f)|^{2}+2\Delta(f)\Big]\bigg \}dvol_{M}\nonumber\\
&&+\int_{\partial M}\Bigg\{\Big[\frac{3}{2}{\rm trace}\Big(\sigma_{n}^{F}-A^{*}(e_{n})\Big)-4h'(0)-{\rm trace}\Big(A(e_{n})\Big)-3{\rm trace}\Big(A^{*}(e_{n})\Big)-2{\rm trace}\Big(\sigma_{n}^{F}\nonumber\\
&&+A(e_{n})\Big)
+(\frac{19i}{16f}+\frac{11}{f}-\frac{15i}{2})\partial_{x_{n}}(f)\Big]\pi {\rm dim}F\Omega_4
+\frac{3\pi g(dx_{n},df)}{8f}\Omega_{4}\Bigg\}dvol_{M}.
\end{eqnarray}
where  $s$ is the scalar curvature.
\end{thm}

\section{Twisted signature operator and its symbol}

Let us recall the definition of twisted signature operators. We consider a $n$-dimensional oriented Riemannian manifold $(M, g^{M})$.
Let $F$ be a real vector bundle over $M$. let $g^{F}$ be an Euclidean metric on $F$. Let
 \begin{equation}
\wedge^{\ast}(T^{\ast}M)=\bigoplus_{i=0}^{n}\wedge^{i}(T^{\ast}M)
\end{equation}
be the real exterior algebra bundle of $T^{\ast}M$. Let
 \begin{equation}
\Omega^{\ast}(M,F)=\bigoplus_{i=0}^{n}\Omega^{i}(M,F)=\bigoplus_{i=0}^{n}C^{\infty}(M,\wedge^{i}(T^{\ast}M)\otimes F)
\end{equation}
be the set of smooth sections of $\wedge^{\ast}(T^{\ast}M)\otimes F$. Let $\ast$ be the Hodge star operator of $g^{TM}$.
It extends  on  $\wedge^{\ast}(T^{\ast}M)\otimes F$ by acting on $F$ as identity. Then $\Omega^{\ast}(M,F)$ inherits the following
standardly induced inner product
 \begin{equation}
\langle \zeta, \eta  \rangle=\int_{M}\langle \zeta\wedge\ast\eta  \rangle_{F},~~~~\zeta, \eta \in\Omega^{\ast}(M,F).
\end{equation}
Let $\widehat{\nabla}^{F}$ be the non-Euclidean connection on $F$. Let $d^{F}$ be the obvious extension of $\nabla^{F}$ on $\Omega^{\ast}(M,F)$.
Let $\delta^{F}=d^{F\ast}$ be the formal adjoint operator of $d^{F}$ with respect to the inner product. Let $\hat{D}^{F}$ be the differential
operator acting on $\Omega^{\ast}(M,F)$ defined by
 \begin{equation}
\hat{D}^{F}=d^{F}+\delta^{F}.
\end{equation}
Let
 \begin{equation}
\omega(F,g^{F})=\widehat{\nabla}^{F,\ast}-\widehat{\nabla}^{F},~~\nabla^{F,e}=\nabla^{F}+\frac{1}{2}\omega(F,g^{F}).
\end{equation}
Then $\nabla^{F,e}$ is an Euclidean connection on $(F,g^{F})$.

Let $\nabla^{\wedge^{\ast}(T^{\ast}M)}$ be the Euclidean connection on $\wedge^{\ast}(T^{\ast}M)$ induced canonically by the Levi-Civita
connection $\nabla^{TM}$ of $g^{TM}$. Let $\nabla^{e}$ be the Euclidean connection on $\wedge^{\ast}(T^{\ast}M)\otimes F$ obtained from
 the tensor product of $\nabla^{\wedge^{\ast}(T^{\ast}M)}$ and $\nabla^{F,e}$.
Let $\{e_{1},\cdots,e_{n}\}$ be an oriented (local) orthonormal basis of $TM$. The following result was proved by Proposition in \cite{BZ}.

The following identity holds
 \begin{equation}
d^{F}+\delta^{F}=\sum_{i=1}^{n}c(e_{i})\nabla^{e}_{e_{i}}-\frac{1}{2}\sum_{i=1}^{n}\hat{c}(e_{i})\omega(F,g^{F})(e_{i}).
\end{equation}

Let $D_{F}^{e}=\sum_{j=1}^{n}c(e_{j})\nabla^{e}_{e_{j}}$ and $\omega(F,g^{F})$ be any element in $\Omega(M,EndF)$, then we define
the generalized twisted signature operators $\hat{D}_{F}$, $\hat{D}^{*}_{F}$ as follows.

For sections $\psi\otimes \chi\in \wedge^{\ast}(T^{\ast}M)\otimes F$,
\begin{eqnarray}\label{eq:1}
&&\hat{D}_{F}(\psi\otimes \chi)=D_{F}^{e}(\psi\otimes \chi)
 -\frac{1}{2}\sum_{i=1}^{n}\hat{c}(e_{i})\omega(F,g^{F})(e_{i})(\psi\otimes \chi),\\
 \label{eq:2}
&&\hat{D}^{*}_{F}(\psi\otimes \chi)=D_{F}^{e}(\psi\otimes \chi)
-\frac{1}{2}\sum_{i=1}^{n}\hat{c}(e_{i})\omega^{*}(F,g^{F})(e_{i})(\psi\otimes \chi).
\end{eqnarray}
Here $\omega^{*}(F,g^{F})(e_{i})$ denotes the adjoint of $\omega(F,g^{F})(e_{i})$.

In the local coordinates $\{x_{i}; 1\leq i\leq n\}$ and the fixed orthonormal frame
$\{\widetilde{e_{1}},\cdots, \widetilde{e_{n}}\}$, the connection matrix $(\omega_{s,t})$ is defined by
\begin{equation}
\widetilde{\nabla}(\widetilde{e_{1}},\cdots, \widetilde{e_{n}})=(\widetilde{e_{1}},\cdots, \widetilde{e_{n}})(\omega_{s,t}).
\end{equation}

 Let $M$ be a $6$-dimensional compact oriented Riemannian manifold with boundary $\partial M$. We define that
$\hat{D}_{F}:~C^{\infty}(M,\wedge^{\ast}(T^{\ast}M)\otimes F)\rightarrow C^{\infty}(M,\wedge^{\ast}(T^{\ast}M)\otimes F)$
is the generalized twisted  signature operator. Take the coordinates and
the orthonormal frame as in Section 3.
 Let $\epsilon (\widetilde{e_j*} ),~\iota (\widetilde{e_j*} )$ be the exterior and interior multiplications respectively. Write
 \begin{equation}
c(\widetilde{e_j})=\epsilon (\widetilde{e_j*} )-\iota (\widetilde{e_j*} );~~
\hat{c}(\widetilde{e_j})=\epsilon (\widetilde{e_j*} )+\iota (\widetilde{e_j*} ).
\end{equation}
  We'll compute ${\rm tr}_{\wedge^*(T^*M)\otimes F}$ in the frame $\{e^{\ast}_{i_1}\wedge\cdots\wedge
e^{\ast}_{i_k}|~1\leq i_1<\cdots<i_k\leq 6\}.$ By (3.2) and (4.8) in \cite{Wa3}, we have
\begin{eqnarray}
 \hat{D}_{F}
    &=&\sum_{i=1}^{n}c(e_{i})\nabla^{e}_{e_{i}}-\frac{1}{2}\sum_{i=1}^{n}\hat{c}(e_{i})\omega(F,g^{F})(e_{i})\nonumber\\
    &=&\sum_{i=1}^{n}c(e_{i})\Big(\nabla_{e_{i}}^{\wedge^{\ast}(T^{\ast}M)}\otimes {\rm id}_{F}+{\rm id}_{\wedge^{\ast}(T^{\ast}M)} \otimes \nabla^{F,e}_{e_{i}} \Big)
    -\frac{1}{2}\sum_{i=1}^{n}\hat{c}(e_{i})\omega(F,g^{F})(e_{i})\nonumber\\
    &=&\sum^n_{i=1}c(\widetilde{e_i})\Big[\widetilde{e_i}+\frac{1}{4}\sum_{s,t}\omega_{s,t}
(\widetilde{e_i})[\hat{c}(\widetilde{e_s})\hat{c}(\widetilde{e_t})-c(\widetilde{e_s})
c(\widetilde{e_t})]\otimes {\rm id}_{F}
  +{\rm id}_{\wedge^{\ast}(T^{\ast}M)}\otimes \sigma^{F,e}_{i}\Big]\nonumber\\
   &&-\frac{1}{2}\sum_{i=1}^{n}\hat{c}(e_{i})\omega(F,g^{F})(e_{i}),
   \end{eqnarray}
  Similarly, we have
\begin{eqnarray}
\hat{D}^{*}_{F}&=&\sum^n_{i=1}c(\widetilde{e_i})\Big[\widetilde{e_i}+\frac{1}{4}\sum_{s,t}\omega_{s,t}
(\widetilde{e_i})[\hat{c}(\widetilde{e_s})\hat{c}(\widetilde{e_t})-c(\widetilde{e_s})c(\widetilde{e_t})]\otimes {\rm id}_{F}+{\rm id}_{\wedge^{\ast}(T^{\ast}M)}\otimes \sigma^{F,e}_{i}\Big]
  \nonumber\\
   &&-\frac{1}{2}\sum_{i=1}^{n}\hat{c}(e_{i})\omega^{*}(F,g^{F})(e_{i}).
\end{eqnarray}
For convenience,
let  $\hat{c}(\omega)=\sum_{i}\hat{c}(e_{i}) \omega(F,g^{F})(e_{i})$
and $\hat{c}(\omega^{*})=\sum_{i}\hat{c}(e_{i}) \omega^{*}(F,g^{F})(e_{i})$, by the composition formula and (2.2.11) in \cite{Wa3}, we obtain in \cite{WJ6},

\begin{lem}
Let $\hat{D}^{*}_{F},  \hat{D}_{F}$ be the twisted signature operators on  $\Gamma(\wedge^*(T^*M)\otimes F)$, then

\begin{eqnarray}
\sigma_1(\hat{D}_{F})&=&\sigma_1(\hat{D}^{*}_{F})=\sqrt{-1}c(\xi);\\
\sigma_0(\hat{D}_{F})&=&\sum^n_{i=1}c(\widetilde{e_i})\Big[\frac{1}{4}\sum_{s,t}\omega_{s,t}
(\widetilde{e_i})[\hat{c}(\widetilde{e_s})\hat{c}(\widetilde{e_t})-c(\widetilde{e_s})
c(\widetilde{e_t})]\otimes {\rm id}_{F}
  +{\rm id}_{\wedge^{\ast}(T^{\ast}M)}\otimes \sigma^{F,e}_{i}\Big]-\frac{\hat{c}(\omega)}{2};\\
\sigma_0(\hat{D}^{*}_{F})&=&\sum^n_{i=1}c(\widetilde{e_i})\Big[\frac{1}{4}\sum_{s,t}\omega_{s,t}
(\widetilde{e_i})[\hat{c}(\widetilde{e_s})\hat{c}(\widetilde{e_t})
-c(\widetilde{e_s})c(\widetilde{e_t})]\otimes {\rm id}_{F}
  +{\rm id}_{\wedge^{\ast}(T^{\ast}M)}\otimes \sigma^{F,e}_{i}\Big]-\frac{\hat{c}(\omega^{*})}{2}.
\end{eqnarray}
\end{lem}

By the composition formula of pseudodifferential operators in Section 2.2.1 of \cite{Wa3}, we have
 \begin{lem}\label{le:31}
The symbol of the  twisted signature operators  $\hat{D}^{*}_{F}, \hat{D}_{F}$ as follows:
\begin{eqnarray}
\sigma_{-1}(\hat{D}_{F}^{-1})&=&\sigma_{-1}((\hat{D}^{*}_{F})^{-1})=\frac{\sqrt{-1}c(\xi)}{|\xi|^{2}}; \\
\sigma_{-2}(\hat{D}_{F}^{-1})&=&\frac{c(\xi)\sigma_{0}(\hat{D}_F)c(\xi)}{|\xi|^{4}}
+\frac{c(\xi)}{|\xi|^{6}}\sum_{j}c(dx_{j})
\Big[\partial_{x_{j}}(c(\xi))|\xi|^{2}-c(\xi)\partial_{x_{j}}(|\xi|^{2})\Big];\\
\sigma_{-2}((\hat{D}^{*}_{F})^{-1})&=&\frac{c(\xi)\sigma_{0}(\hat{D}^{*}_F)c(\xi)}{|\xi|^{4}}
+\frac{c(\xi)}{|\xi|^{6}}\sum_{j}c(dx_{j})
\Big[\partial_{x_{j}}(c(\xi))|\xi|^{2}-c(\xi)\partial_{x_{j}}(|\xi|^{2})\Big].
\end{eqnarray}
\end{lem}

Since $\Psi$ is a global form on $\partial M$, so for any fixed point $x_{0}\in\partial M$, we can choose the normal coordinates
$U$ of $x_{0}$ in $\partial M$(not in $M$) and compute $\Psi(x_{0})$ in the coordinates $\widetilde{U}=U\times [0,1)$ and the metric
$\frac{1}{h(x_{n})}g^{\partial M}+dx _{n}^{2}$. The dual metric of $g^{\partial M}$ on $\widetilde{U}$ is
$\frac{1}{\tilde{h}(x_{n})}g^{\partial M}+dx _{n}^{2}.$ Write
$g_{ij}^{M}=g^{M}(\frac{\partial}{\partial x_{i}},\frac{\partial}{\partial x_{j}})$;
$g^{ij}_{M}=g^{M}(d x_{i},dx_{j})$, then

\begin{equation}
[g_{i,j}^{M}]=
\begin{bmatrix}\frac{1}{h( x_{n})}[g_{i,j}^{\partial M}]&0\\0&1\end{bmatrix};\quad
[g^{i,j}_{M}]=\begin{bmatrix} h( x_{n})[g^{i,j}_{\partial M}]&0\\0&1\end{bmatrix},
\end{equation}
and
\begin{equation}
\partial_{x_{s}} g_{ij}^{\partial M}(x_{0})=0,\quad 1\leq i,j\leq n-1;\quad g_{i,j}^{M}(x_{0})=\delta_{ij}.
\end{equation}

Let $\{e_{1},\cdots, e_{n-1}\}$ be an orthonormal frame field in $U$ about $g^{\partial M}$ which is parallel along geodesics and
$e_{i}=\frac{\partial}{\partial x_{i}}(x_{0})$, then $\{\widetilde{e_{1}}=\sqrt{h(x_{n})}e_{1}, \cdots,
\widetilde{e_{n-1}}=\sqrt{h(x_{n})}e_{n-1},\widetilde{e_{n}}=\texttt{d}x_{n}\}$ is the orthonormal frame field in $\widetilde{U}$ about $g^{M}.$
Locally $\wedge^{\ast}(T^{\ast}M)|\widetilde{U}\cong \widetilde{U}\times\wedge^{*}_{C}(\frac{n}{2}).$ Let $\{f_{1},\cdots,f_{n}\}$ be the orthonormal basis of
$\wedge^{*}_{C}(\frac{n}{2})$. Take a spin frame field $\sigma: \widetilde{U}\rightarrow Spin(M)$ such that
$\pi\sigma=\{\widetilde{e_{1}},\cdots, \widetilde{e_{n}}\}$ where $\pi: Spin(M)\rightarrow O(M)$ is a double covering, then
$\{[\sigma, f_{i}], 1\leq i\leq 4\}$ is an orthonormal frame of $\wedge^{\ast}(T^{\ast}M)|_{\widetilde{U}}.$ In the following,
since the global form $\Psi$
is independent of the choice of the local frame, so we can compute $\rm{trace}_{\wedge^{\ast}(T^{\ast}M)}$ in the frame
$\{[\sigma, f_{i}], 1\leq i\leq 4\}$.
Let $\{E_{1},\cdots,E_{n}\}$ be the canonical basis of $R^{n}$ and
$c(E_{i})\in cl_{C}(n)\cong \rm{Hom}(\wedge^{*}_{C}(\frac{n}{2}),\wedge^{*}_{C}(\frac{n}{2}))$ be the Clifford action. By \cite{Wa3}, then

\begin{equation}
c(\widetilde{e_{i}})=\bigg[\big(\sigma,c(E_{i})\big)\bigg]; \quad c(\widetilde{e_{i}})[(\sigma, f_{i})]=[\sigma,(c(E_{i}))f_{i}]; \quad
\frac{\partial}{\partial x_{i}}=[(\sigma,\frac{\partial}{\partial x_{i}})],
\end{equation}
then we have $\frac{\partial}{\partial x_{i}}c(\widetilde{e_{i}})=0$ in the above frame. By Lemma 2.2 in \cite{Wa3}, we have

\begin{lem}\label{le:32}
\begin{eqnarray}
\partial_{x_j}(|\xi|_{g^M}^2)(x_0)&=&\left\{
       \begin{array}{c}
        0,  ~~~~~~~~~~ ~~~~~~~~~~ ~~~~~~~~~~~~~{\rm if }~j<n; \\[2pt]
       h'(0)|\xi'|^{2}_{g^{\partial M}},~~~~~~~~~~~~~~~~~~~~~{\rm if }~j=n,
       \end{array}
    \right. \\
\partial_{x_j}[c(\xi)](x_0)&=&\left\{
       \begin{array}{c}
      0,  ~~~~~~~~~~ ~~~~~~~~~~ ~~~~~~~~~~~~~{\rm if }~j<n;\\[2pt]
\partial _{x_{n}}(c(\xi'))(x_{0}), ~~~~~~~~~~~~~~~~~{\rm if }~j=n,
       \end{array}
    \right.
\end{eqnarray}
where $\xi=\xi'+\xi_{n}dx_{n}$.
\end{lem}
Then an application of Lemma 2.3 in \cite{Wa3} shows
\begin{lem}
The symbol of the  twisted signature operators  $\hat{D}^{*}_{F},   \hat{D}_{F}$ as follows:
\begin{eqnarray}
\sigma_{0}(\hat{D}^{*}_{F})&=&\theta+\vartheta^*;\\
\sigma_{0}(\hat{D}_{F})&=&\theta+\vartheta,
\end{eqnarray}
where
\begin{eqnarray}
\theta&=&-\frac{5}{4}h'(0)c(dx_n)
+\frac{1}{4}h'(0)\sum^{n-1}_{i=1}c(\widetilde{e_i})\hat{c}(\widetilde{e_n})\hat{c}(\widetilde{e_i})(x_0)\otimes {\rm id}_{F};\nonumber\\
\vartheta^*&=&\sum^n_{i=1}c(\widetilde{e_i})\sigma^{F,e}_{i}-\frac{1}{2}
\sum_{i=1}^{n}\hat{c}(e_{i})\omega^{*}(F,g^{F})(e_{i});\nonumber\\
\vartheta&=&\sum^n_{i=1}c(\widetilde{e_i})\sigma^{F,e}_{i}-\frac{1}{2}\sum_{i=1}^{n}\hat{c}(e_{i})\omega(F,g^{F})(e_{i}).
\end{eqnarray}
\end{lem}

In order to get the symbol of operators $\hat{D}^{*}_{F}f\cdot
\hat{D}_{F}f^{-1}\cdot\hat{D}_{F}^{*}f$. Similar to (3.19)-(3.23), we give the specification of
$\hat{D}^{*}_{F}f\cdot\hat{D}_{F}f^{-1}\cdot\hat{D}_{F}^{*}f$.

Combining (4.11) and (4.12), we have
\begin{eqnarray*}
&&\hat{D}^{*}_{F}f\cdot\hat{D}_{F}f^{-1}\cdot\hat{D}_{F}^{*}f\nonumber\\
&=&f\cdot\hat{D}_{F}^{*}\hat{D}_{F}\hat{D}_{F}^{*}
+c(df)\hat{D}_{F}\hat{D}_{F}^{*}
-\hat{D}_{F}^{*}\hat{D}_{F}f\cdot c(df^{-1})\cdot f
+\hat{D}_{F}^{*}\cdot c(df) c(df^{-1})f\nonumber\\
&=&f\cdot\Bigg\{\sum^{n}_{i,j,l=1}\sum^{n}_{r=1}c(e_{r})\langle e_{r},dx_{l}\rangle(-g^{ij}\partial_{l}\partial_{i}\partial_{j})
         +\sum^{n}_{r,l=1}c(e_{r})\langle e_{r},dx_{l}\rangle \bigg\{ -\sum^{n}_{i,j=1}(\partial_{l}g^{ij})\partial_{i}\partial_{j}-
         \sum^{n}_{i,k,j=1}g^{ij}\nonumber\\
         &&\times(4\sigma^{\wedge^*(T^*M)\otimes F}_{i}\partial_{j}-2\Gamma^{k}_{ij} \times\partial_{k})\partial_{l}\bigg\}+\sigma_{0}(\hat{D}^{*}_{F})
         (-\sum^{n}_{i,j=1}g^{ij}\partial_{i}\partial_{j})-\sum^{n}_{r,l=1}c(e_{r})\langle e_{r},dx_{l}\rangle \nonumber\\
         &&\times\sum^{n}_{j,k=1}\Big[\hat{c}(w)c(e_{j})+c(e_{j})\hat{c}(w^{*})\Big]\langle e_{j},dx^{k}\rangle
         \times\partial_{l}\partial_{k}+\sum^{n}_{r,l=1}c(e_{r})\langle e_{r},dx_{l}\rangle\partial_{l}  \bigg\{-\sum^{n}_{i,j,k=1}g^{ij} \nonumber\\
         && \times\Big[(\partial_{i}
         \sigma^{j,e}_{\wedge^{\ast}(T^{\ast}M)\otimes F})+\sigma^{i}_{\wedge^{\ast}(T^{\ast}M)\otimes F}
         \sigma^{j,e}_{\wedge^{\ast}(T^{\ast}M)\otimes F,e}-\Gamma_{ij}^{k}\sigma_{\wedge^{\ast}(T^{\ast}M)\otimes F}^{k}\Big]+\frac{1}{4}s-\frac{1}{2}
         \sum^{n}_{j=1}\hat{c}(\omega)c(e_{j})\nonumber\\
        &&\times\sigma^{\wedge^{\ast}(T^{\ast}M)\otimes F,e}_{j}
        -\frac{1}{2}\sum^{n}_{j=1}c(e_{j})e_{j}\big(\hat{c}(\omega^{*})\big)+\frac{1}{2}\sum_{i\neq j} R^{F,e}(e_{i},e_{j})c(e_{i})c(e_{j})+\frac{1}{4}\hat{c}(\omega)\hat{c}(\omega^{*}) -\frac{1}{2}\sum^{n}_{j=1}c(e_{j})\nonumber\\
        &&\sigma_{j}^{\wedge^{\ast}(T^{\ast}M)\otimes F,e} \hat{c}(\omega^{*})
        \bigg\}+\sigma_{0}(\hat{D}^{*}_{F})\bigg\{-2\sigma^{j}_{\wedge^{\ast}(T^{\ast}M)\otimes F}\partial_{j}+\Gamma^{k}\partial_{k}
        -\frac{1}{2}\sum^{n}_{j=1}\Big[\hat{c}(\omega)c(e_{j})+c(e_{j})\hat{c}(\omega^{*}) \Big]\nonumber\\
        &&\times e_{j}-g^{ij}\Big[(\partial_{i}\sigma^{j,e}_{\wedge^{\ast}(T^{\ast}M)\otimes F}) +\sigma^{i}_{\wedge^{\ast}(T^{\ast}M)\otimes F}
        \sigma^{j,e}_{\wedge^{\ast}(T^{\ast}M)\otimes F,e}
        -\Gamma_{ij}^{k}\sigma_{\wedge^{\ast}(T^{\ast}M)\otimes F}^{k}\Big]+\frac{1}{4}\hat{c}(\omega)\hat{c}(\omega^{*}) \nonumber\\
         &&-\frac{1}{2}\sum^{n}_{j=1}\hat{c}(\omega)c(e_{j})\sigma^{\wedge^{\ast}(T^{\ast}M)\otimes F,e}_{j}
          -\frac{1}{2}\sum^{n}_{j=1}c(e_{j})e_{j}\big(\hat{c}(\omega^{*})\big)
           -\frac{1}{2}\sum^{n}_{j=1}c(e_{j})\sigma_{j}^{\wedge^{\ast}(T^{\ast}M)\otimes F,e} \hat{c}(\omega^{*})+\frac{1}{4}s \nonumber\\
          &&+\frac{1}{2}\sum_{i\neq j} R^{F,e}(e_{i},e_{j})c(e_{i})c(e_{j})\bigg\}
          +\sum^{n}_{r,l=1}c(e_{r})\langle e_{r},dx_{l}\rangle
          \bigg\{\sum^{n}_{i,j,k=1}g^{ij}(\partial_{l}\Gamma^{k}_{ij})\partial_{k}
          -2\sum^{n}_{i,j=1}g^{ij}(\partial_{l}\sigma^{\wedge^*(T^*M)\otimes F}_{i})\nonumber\\
          &&\times\partial_{j}-2\sum^{n}_{i,j=1}(\partial_{l}g^{ij})\sigma^{\wedge^*(T^*M)\otimes F}_{i}\partial_{j}-\frac{1}{2}\sum^{n}_{j,k=1}\Big[\partial_{l}
          \Big(\hat{c}(w)c(e_{j})+c(e_{j})\hat{c}(w^{*})\Big)\Big]\langle e_{j},dx^{k}\rangle\partial_{k}\nonumber\\
  \end{eqnarray*}
\begin{eqnarray}
          &&+\sum^{n}_{i,j,k=1}(\partial_{l}g^{ij})\Gamma^{k}_{ij}\partial_{k}
          -\frac{1}{2}\sum^{n}_{j,k=1}\Big(\hat{c}(w)c(e_{j})
          +c(e_{j})\hat{c}(w^{*})\Big)\Big[\partial_{l}\langle e_{j},dx^{k}\rangle\Big]\partial_{k} \bigg\}\Bigg\}+c(df)\Bigg\{-g^{ij}\partial_{i}\partial_{j}\nonumber\\
          &&-2\sigma^{j}_{\wedge^{\ast}(T^{\ast}M)\otimes F}\partial_{j}+\Gamma^{k}\partial_{k}
                -\frac{1}{2}\sum_{j}\Big[\hat{c}(\omega)c(e_{j})+c(e_{j})\hat{c}(\omega^{*}) \Big]e_{j}-g^{ij}\Big[(\partial_{i}\sigma^{j,e}_{\wedge^{\ast}(T^{\ast}M)\otimes F})+\sigma^{i}_{\wedge^{\ast}(T^{\ast}M)\otimes F}\nonumber\\
                &&
                \sigma^{j,e}_{\wedge^{\ast}(T^{\ast}M)\otimes F,e}
                  -\Gamma_{ij}^{k}\sigma_{\wedge^{\ast}(T^{\ast}M)\otimes F}^{k}\Big]-\frac{1}{2}\sum_{j}\hat{c}(\omega)c(e_{j})\sigma^{\wedge^{\ast}(T^{\ast}M)\otimes F,e}_{j}
               -\frac{1}{2}\sum_{j}c(e_{j})e_{j}\big(\hat{c}(\omega^{*})\big)+\frac{1}{4}s\nonumber\\
                &&+\frac{1}{4}\hat{c}(\omega)\hat{c}(\omega^{*}) -\frac{1}{2}\sum_{j}c(e_{j})\sigma_{j}^{\wedge^{\ast}(T^{\ast}M)\otimes F,e} \hat{c}(\omega^{*})
                 +\frac{1}{2}\sum_{i\neq j} R^{F,e}(e_{i},e_{j})c(e_{i})c(e_{j})\Bigg\}-\Bigg\{-g^{ij}\partial_{i}\partial_{j}\nonumber\\
&&-2\sigma^{j}_{\wedge^{\ast}(T^{\ast}M)\otimes F}\partial_{j}+\Gamma^{k}\partial_{k}
                -\frac{1}{2}\sum_{j}\Big[\hat{c}(\omega)c(e_{j})+c(e_{j})\hat{c}(\omega^{*}) \Big]e_{j}-g^{ij}\Big[(\partial_{i}\sigma^{j,e}_{\wedge^{\ast}(T^{\ast}M)\otimes F}) +\sigma^{i}_{\wedge^{\ast}(T^{\ast}M)\otimes F}\nonumber\\
                &&\times
                \sigma^{j,e}_{\wedge^{\ast}(T^{\ast}M)\otimes F,e}
                  -\Gamma_{ij}^{k}\sigma_{\wedge^{\ast}(T^{\ast}M)\otimes F}^{k}\Big]
                  -\frac{1}{2}\sum_{j}\hat{c}(\omega)c(e_{j})\sigma^{\wedge^{\ast}(T^{\ast}M)\otimes F,e}_{j}-\frac{1}{2}\sum_{j}c(e_{j})e_{j}\big(\hat{c}(\omega^{*})\big)\nonumber\\
                &&
+\frac{1}{4}s+\frac{1}{4}\hat{c}(\omega)\hat{c}(\omega^{*}) -\frac{1}{2}\sum_{j}c(e_{j})\sigma_{j}^{\wedge^{\ast}(T^{\ast}M)\otimes F,e} \hat{c}(\omega^{*})
                 +\frac{1}{2}\sum_{i\neq j} R^{F,e}(e_{i},e_{j})c(e_{i})c(e_{j})\Bigg\}
f\cdot c(df^{-1})\cdot f \nonumber\\
&&+\Bigg\{\sum^{n}_{i,j=1}g^{ij}c(\partial_{i})
\Big[\partial_{j}+\Big(\frac{1}{4}\sum_{s,t}\omega_{s,t}
(\widetilde{e_i})[\hat{c}(\widetilde{e_s})\hat{c}(\widetilde{e_t})
-c(\widetilde{e_s})c(\widetilde{e_t})]\otimes id_{F}
  +id_{\wedge^{\ast}(T^{\ast}M)}\otimes \sigma^{F,e}_{i}\Big)\Big]\nonumber\\
  &&-\frac{\hat{c}(\omega^{*})}{2}\Bigg\}\cdot c(df) c(df^{-1})f.
\end{eqnarray}

By the above composition formulas, then we obtain:

\begin{lem}
Let $\hat{D}^{*}_{F}, \hat{D}_{F}$ be the twisted signature operators on  $\Gamma(\wedge^*(T^*M)\otimes F)$, then
\begin{eqnarray}
&&\sigma_{3}(\hat{D}^{*}_{F}f\cdot\hat{D}_{F}f^{-1}\cdot\hat{D}^{*}_{F}f)
=f\sigma_{3}(\hat{D}^{*}_{F}\hat{D}_{F}\hat{D}^{*}_{F})=f\sqrt{-1}c(\xi)|\xi|^2; \\
&&\sigma_{2}(\hat{D}^{*}_{F}f\cdot\hat{D}_{F}f^{-1}\cdot\hat{D}^{*}_{F}f)
=f\sigma_{2}(\hat{D}^{*}_{F}\hat{D}_{F}\hat{D}^{*}_{F})+2c(df)|\xi|^2
\end{eqnarray}
where $\sigma_{2}(\hat{D}^{*}_{F}\hat{D}_{F}\hat{D}^{*}_{F})
=c(\xi)(4\sigma^k-2\Gamma^k)\xi_{k}
     -\frac{1}{4}|\xi|^2h'(0)c(dx_{n})
     +|\xi|^2\Big(\frac{1}{4}h'(0)\sum\limits^{5}_{i=1}c(\widetilde{e_i})
     \hat{c}(\widetilde{e_n})\hat{c}(\widetilde{e_i})(x_0)+\vartheta^*-\hat{c}(w^*)\Big)
     +c(\xi)\hat{c}(w)c(\xi)$.
\end{lem}
For convenience,
we write that
$\sigma_{2}(\hat{D}^{*}_{F}\hat{D}_{F}\hat{D}^{*}_{F})
=G+|\xi|^2\Big(p+\vartheta^*-\hat{c}(w^*)\Big)
     +c(\xi)\hat{c}(w)c(\xi)$.
By (4.28), (4.29), Lemma 2.1 in \cite{Wa3} and  the composition formula of psudodifferential operators, similar to (3.26)-(3.28), we obtain

\begin{lem}Let $\hat{D}^{*}_{F}, \hat{D}_{F}$ be the generalized twisted signature operators on  $\Gamma(\wedge^{\ast}(T^{\ast}M)\otimes F)$, then
\begin{eqnarray}
\sigma_{-3}(\hat{D}^{*}_{F}f\cdot\hat{D}_{F}f^{-1}\cdot\hat{D}^{*}_{F}f)^{-1}
&=&\frac{\sqrt{-1}c(\xi)}{f|\xi|^4}; \\
\sigma_{-4}(\hat{D}^{*}_{F}f\cdot\hat{D}_{F}f^{-1}\cdot\hat{D}^{*}_{F}f)^{-1}
&=&f^{-1}\sigma_{-4}\big((\hat{D}^{*}_{F}\hat{D}_{F}\hat{D}^{*}_{F})^{-1}\big)
+\frac{2c(\xi)c(df)c(\xi)}{f^{2}|\xi|^6}\nonumber\\
&&+\frac{ic(\xi)\sum\limits_j\Big[c(dx_j)|\xi|^2+2\xi_{j}c(\xi)\Big]D_{x_j}(f^{-1})c(\xi)}{|\xi|^8},
\end{eqnarray}
where
\begin{eqnarray}
&&\sigma_{-4}\big((\hat{D}^{*}_{F}\hat{D}_{F}\hat{D}^{*}_{F})^{-1}\big)\nonumber\\
&=&\frac{c(\xi)\sigma_{2}(\hat{D}^{*}_{F}\hat{D}_{F}\hat{D}^{*}_{F})c(\xi)}{|\xi|^8}
+\frac{c(\xi)}{|\xi|^{10}}\sum_j\Big[c(dx_j)|\xi|^2+2\xi_{j}c(\xi)\Big]
\Big[\partial_{x_j}[c(\xi)]|\xi|^2-2c(\xi)\partial_{x_j}(|\xi|^2)\Big]\nonumber\\
&=&\frac{c(\xi)Gc(\xi)}{|\xi|^8}+\frac{c(\xi)\Big(p+\vartheta^*-\hat{c}(w^*)\Big)c(\xi)}{|\xi|^6}
+\frac{\hat{c}(w)}{|\xi|^4}+\frac{c(\xi)}{|\xi|^{10}}\sum_j\Big[c(dx_j)|\xi|^2+2\xi_{j}c(\xi)\Big]
\Big[\partial_{x_j}[c(\xi)]|\xi|^2\nonumber\\
&&-2c(\xi)\partial_{x_j}(|\xi|^2)\Big]
.
\end{eqnarray}
\end{lem}

Hence we cite that
\begin{thm} \cite{WJ6}
For even $n$-dimensional oriented  compact Riemainnian manifolds without boundary,
the following equality holds:
\begin{eqnarray}
&&Wres(\hat{D}^{*}_{F}f\cdot\hat{D}_{F}f^{-1})^{(\frac{-n+2}{2})}\nonumber\\
&=&\frac{(2\pi)^{\frac{n}{2}}}{(\frac{n}{2}-2)!}\int_{M}\bigg\{{\rm{trace}}
 \Big[-\frac{s}{12} +\frac{n}{16}\big[\hat{c}(\omega^{*})-\hat{c}(\omega)\big]^{2}
         -\frac{1}{4}\hat{c}(\omega^{*})\hat{c}(\omega) -\frac{1}{4}\sum_{j}\nabla_{e_{j}}^{F}\big(\hat{c}(\omega^{*})\big)c(e_{j})\nonumber\\
       &&+\frac{1}{4}\sum_{j}c(e_{j})\nabla_{e_{j}}^{F}\big(\hat{c}(\omega)\big) \Big]+4f^{-1}\Delta(f)+8\big\langle grad_{M}f,grad_{M}(f^{-1})\big\rangle-5f^{-2}\Big[|grad_{M}f|^2\nonumber\\
       &&+2\Delta f\Big]\bigg\}dvol_{M}.
\end{eqnarray}
\end{thm}

 \section{Conformal perturbations of twisted Signature Operators and Noncommutative residue}
In the following, we will compute the more general case $\widetilde{Wres}[\pi^{+}(f\hat{D}_{F}^{-1}) \circ\pi^{+}\big(f^{-1}(\hat{D}_{F}^{*})^{-1}\cdot f\hat{D}_{F}^{-1}\cdot f^{-1}(\hat{D}^{*}_{F})^{-1}\big)]$ for nonzero
smooth functions $f,~f^{-1}$.
An application of (2.1.4) in \cite{Wa5} shows that

\begin{eqnarray}
&&\widetilde{Wres}\big[\pi^{+}(f\hat{D}_{F}^{-1}) \circ\pi^{+}\big(f^{-1}(\hat{D}_{F}^{*})^{-1}\cdot f\hat{D}_{F}^{-1}\cdot f^{-1}(\hat{D}^{*}_{F})^{-1}\big)\big]\nonumber\\
&=&\int_{M}\int_{|\xi|=1}\rm{trace}_{\wedge^{\ast}(T^{\ast}M)\otimes F}
  \big((\hat{\emph{D}}^{*}_{\emph{F}}\emph{f}\cdot\hat{\emph{D}}_{\emph{F}}\emph{f}^{-1})^{-2}\big)
  \sigma(\xi)\emph{dx}+\int_{\partial M} \Psi,
\end{eqnarray}
where
 \begin{eqnarray}
 \Psi&=&\int_{|\xi'|=1}\int_{-\infty}^{+\infty}\sum_{j,k=0}^{\infty}\sum \frac{(-i)^{|\alpha|+j+k+\ell}}{\alpha!(j+k+1)!}
\rm{trace}_{\wedge^{\ast}(T^{\ast}M)\otimes F}\Big[\partial_{x_{n}}^{j}\partial_{\xi'}^{\alpha}\partial_{\xi_{n}}^{k}
\sigma_{r}^{+}(\emph{f}~\hat{\emph{D}}_{\emph{F}}^{-1})(x',0,\xi',\xi_{n})\nonumber\\
&&\times\partial_{x_{n}}^{\alpha}\partial_{\xi_{n}}^{j+1}\partial_{x_{n}}^{k}
\sigma_{l}\big(f^{-1}(\hat{D}_{F}^{*})^{-1}\cdot f\hat{D}_{F}^{-1}\cdot f^{-1}(\hat{D}^{*}_{F})^{-1}\big)(x',0,\xi',\xi_{n})\Big]
d\xi_{n}\sigma(\xi')dx' ,
\end{eqnarray}
and the sum is taken over $r-k+|\alpha|+\ell-j-1=-n,r\leq-1,\ell\leq-1$.

Locally we can use Theorem 4.7 to compute the interior term of (5.1), then
\begin{eqnarray}
&&\int_{M}\int_{|\xi|=1}\rm{trace}_{\wedge^{\ast}(T^{\ast}M)\otimes F}
  [\sigma_{-4}((\hat{\emph{D}}_{\emph{F}}^{*}\emph{f}\cdot\hat{\emph{D}}_{\emph{F}}\emph{f}^{-1})^{-2})]\sigma(\xi)
  \emph{dx}\nonumber\\
&=&8\pi^{3}\int_{M}\Bigg\{\rm{trace}
 \Big[-\frac{s}{12} +\frac{3}{8}\big[\hat{c}(\omega^{*})-\hat{c}(\omega)\big]^{2}
         -\frac{1}{4}\hat{c}(\omega^{*})\hat{c}(\omega)
         -\frac{1}{4}\sum_{j}\nabla_{e_{j}}^{F}\big(\hat{c}(\omega^{*})\big)c(e_{j})\nonumber\\
       &&  +\frac{1}{4}\sum_{j}c(e_{j})\nabla_{e_{j}}^{F}\big(\hat{c}(\omega)\big)\Big]
       +4f^{-1}\Delta(f)
       +8\big\langle grad_M(f),grad_M(f^{-1})\big\rangle-5f^{-2}\Big[|grad_M(f)|^2\nonumber\\
       &&+2\Delta(f)\Big]\Bigg\}dvol_{M}.
\end{eqnarray}

So we only need to compute $\int_{\partial M} \Psi$. From the remark above, now we can compute $\Psi$ (see formula (5.2) for the definition of $\Psi$).
Since the sum is taken over $r+\ell-k-j-|\alpha|-1=-6, \ r\leq-1, \ell\leq -3$, then we have the $\int_{\partial_{M}}\Psi$
is the sum of the following five cases:
~\\
~\\
\noindent  {\bf case (a)~(I)}~$r=-1, l=-3, j=k=0, |\alpha|=1$.\\

By (5.2), we get
 \begin{eqnarray}
{\rm case~(a)~(I)}&=&-\int_{|\xi'|=1}\int^{+\infty}_{-\infty}\sum_{|\alpha|=1}{\rm trace}
\Big[\partial^{\alpha}_{\xi'}\pi^{+}_{\xi_{n}}\sigma_{-1}(f\hat{D}_{F}^{-1})
      \times\partial^{\alpha}_{x'}\partial_{\xi_{n}}\sigma_{-3}
      \big(f^{-1}(\hat{D}_{F}^{*})^{-1}\cdot f\hat{D}_{F}^{-1}\nonumber\\
      &&\cdot f^{-1}(\hat{D}^{*}_{F})^{-1}\big)\Big](x_0)d\xi_n\sigma(\xi')dx'\nonumber\\
 &=&-\int_{|\xi'|=1}\int^{+\infty}_{-\infty}\sum_{|\alpha|=1}{\rm trace}
\Big[\partial^{\alpha}_{\xi'}\pi^{+}_{\xi_{n}}\sigma_{-1}(\hat{D}_{F}^{-1})
      \times\partial^{\alpha}_{x'}\partial_{\xi_{n}}\sigma_{-3}
(\hat{D}_{F}^{*}\hat{D}_{F}\hat{D}^{*}_{F})^{-1}\Big](x_0)d\xi_n\sigma(\xi')dx'\nonumber\\
      &&-f\sum\limits_{j<n}\partial_{j}(f^{-1})\int_{|\xi'|=1}\int^{+\infty}_{-\infty}\sum_{|\alpha|=1}{\rm trace}
\Big[\partial^{\alpha}_{\xi'}\pi^{+}_{\xi_{n}}\sigma_{-1}(\hat{D}_{F}^{-1})
      \times\partial_{\xi_{n}}\sigma_{-3}
(\hat{D}_{F}^{*}\hat{D}_{F}\hat{D}^{*}_{F})^{-1}\Big](x_0)\nonumber\\
      &&\times d\xi_n\sigma(\xi')dx'.
\end{eqnarray}
By (3.24) and (4.29), we have
$\sigma_{-3}\big((\hat{D}^{*}_{F}\hat{D}_{F}\hat{D}^{*}_{F})^{-1}\big)
=\sigma_{-3}\big((\widetilde{D}^{*}_{F}\widetilde{D}_{F}\widetilde{D}^{*}_{F})^{-1}\big)$.

By (3.34) and Lemma 2.2 in \cite{Wa3}, for $i<n$ we have
 \begin{equation}
 \partial_{x_{i}}\sigma_{-3}\big((\hat{D}^{*}_{F}\hat{D}_{F}\hat{D}^{*}_{F})^{-1}\big)(x_0)=0.
\end{equation}
 Thus we have
\begin{eqnarray}
-\int_{|\xi'|=1}\int^{+\infty}_{-\infty}\sum_{|\alpha|=1}{\rm trace}
\Big[\partial^{\alpha}_{\xi'}\pi^{+}_{\xi_{n}}\sigma_{-1}(\hat{D}_{F}^{-1})
      \times\partial^{\alpha}_{x'}\partial_{\xi_{n}}\sigma_{-3}
(\hat{D}_{F}^{*}\hat{D}_{F}\hat{D}^{*}_{F})^{-1}\Big](x_0)d\xi_n\sigma(\xi')dx'=0.
\end{eqnarray}
By (3.12) and (4.16), we have
$\sigma_{-1}(\hat{D}_{F})^{-1}
=\sigma_{-1}(\widetilde{D}_{F})^{-1}$.
Similar to (3.36)-(3.38), for $i<n$, we have
\begin{eqnarray}
&&{\rm trace}
\Big[\partial^{\alpha}_{\xi'}\pi^{+}_{\xi_{n}}\sigma_{-1}(\hat{D}_{F}^{-1})
      \times\partial_{\xi_{n}}\sigma_{-3}
(\hat{D}_{F}^{*}\hat{D}_{F}\hat{D}^{*}_{F})^{-1}\Big](x_0)\nonumber\\
&=&-\xi_i{\rm trace}
\Big[\frac{c(dx_n)^{2}}{2(\xi_n-\sqrt{-1})^2}\Big]-4\sqrt{-1}\xi_n\xi_i{\rm trace}
\Big[\frac{c(dx_i)^{2}}{2(\xi_n-\sqrt{-1})|\xi|^6}\Big]+4\sqrt{-1}\xi_n\xi_i(\xi_n\nonumber\\
&&-2\sqrt{-1}){\rm trace}
\Big[\frac{c(\xi')^{2}}{2(\xi_n-\sqrt{-1})^2|\xi|^6}\Big]+4\sqrt{-1}\xi^{2}_n\xi_i{\rm trace}
\Big[\frac{c(dx_n)^{2}}{2(\xi_n-\sqrt{-1})^2|\xi|^6}\Big].
\end{eqnarray}
We note that $i<n,~\int_{|\xi'|=1}\xi_i\sigma(\xi')=0$,
so
\begin{eqnarray}
&&-f\sum\limits_{j<n}\partial_{j}(f^{-1})\int_{|\xi'|=1}\int^{+\infty}_{-\infty}\sum_{|\alpha|=1}{\rm trace}
\Big[\partial^{\alpha}_{\xi'}\pi^{+}_{\xi_{n}}\sigma_{-1}(\hat{D}_{F}^{-1})
      \times\partial_{\xi_{n}}\sigma_{-3}
(\hat{D}_{F}^{*}\hat{D}_{F}\hat{D}^{*}_{F})^{-1}\Big](x_0) d\xi_n\sigma(\xi')dx'\nonumber\\
&=&0.
\end{eqnarray}
Then we have ${\bf case~(a)~(I)}=0$.
~\\

\noindent  {\bf case (a)~(II)}~$r=-1, l=-3, |\alpha|=k=0, j=1$.\\

By (5.2), we have
 \begin{eqnarray}
{\rm case~(a)~(II)}&=&-\frac{1}{2}\int_{|\xi'|=1}\int^{+\infty}_{-\infty} {\rm
trace} \Big[\partial_{x_{n}}\pi^{+}_{\xi_{n}}\sigma_{-1}(f\hat{D}_{F}^{-1})
\times\partial^{2}_{\xi_{n}}\sigma_{-3}
\big(f^{-1}(\hat{D}_{F}^{*})^{-1}\cdot f\hat{D}_{F}^{-1}\cdot\nonumber\\
&&f^{-1}(\hat{D}^{*}_{F})^{-1}\big)\Big](x_0)d\xi_n\sigma(\xi')dx'\nonumber\\
  &=&-\frac{1}{2}\int_{|\xi'|=1}\int^{+\infty}_{-\infty}{\rm trace}
\Big[\partial_{x_{n}}\pi^{+}_{\xi_{n}}\sigma_{-1}(\hat{D}_{F}^{-1})
      \times\partial^2_{\xi_{n}}\sigma_{-3}
(\hat{D}_{F}^{*}\hat{D}_{F}\hat{D}^{*}_{F})^{-1}\Big](x_0)d\xi_n\sigma(\xi')dx'\nonumber\\
      &&-\frac{1}{2}f^{-1}\partial_{x_{n}}(f)
      \int_{|\xi'|=1}\int^{+\infty}_{-\infty}{\rm trace}
\Big[\pi^{+}_{\xi_{n}}\sigma_{-1}(\hat{D}_{F}^{-1})
      \times\partial^2_{\xi_{n}}\sigma_{-3}
(\hat{D}_{F}^{*}\hat{D}_{F}\hat{D}^{*}_{F})^{-1}\Big](x_0)\nonumber\\
      &&\times d\xi_n\sigma(\xi')dx'.
\end{eqnarray}
Since $n=6$, ${\rm trace}_{\wedge^*(T^*M)}[-{\rm id}]=-64{\rm dim}F$. By the relation of the Clifford action and ${\rm trace}PQ={\rm trace}QP$,  then
\begin{eqnarray}
&&{\rm trace}[c(\xi')c(dx_{n})]=0; \ {\rm trace}[c(dx_{n})^{2}]=-64{\rm dim}F;\
{\rm trace}[c(\xi')^{2}](x_{0})|_{|\xi'|=1}=-64{\rm dim}F;\nonumber\\
&&{\rm trace}[\partial_{x_{n}}[c(\xi')]c(dx_{n})]=0; \
{\rm trace}[\partial_{x_{n}}c(\xi')c(\xi')](x_{0})|_{|\xi'|=1}=-32h'(0){\rm dim}F.
\end{eqnarray}

Similar to (3.41)-(3.45), then we obtain

\begin{eqnarray}
&&-\frac{1}{2}\int_{|\xi'|=1}\int^{+\infty}_{-\infty}{\rm trace}
\Big[\partial_{x_{n}}\pi^{+}_{\xi_{n}}\sigma_{-1}(\hat{D}_{F}^{-1})
      \times\partial^2_{\xi_{n}}\sigma_{-3}
(\hat{D}_{F}^{*}\hat{D}_{F}\hat{D}^{*}_{F})^{-1}\Big](x_0)d\xi_n\sigma(\xi')dx'\nonumber\\
&=&-\frac{1}{2}\int_{|\xi'|=1}\int^{+\infty}_{-\infty}
8h'(0){\rm dim}F\frac{-8-24\xi_{n}i+40\xi^{2}_{n}+24i\xi^{3}_{n}}
{(\xi_{n}-i)^{6}(\xi_{n}+i)^{4}}d\xi_n\sigma(\xi')dx'\nonumber\\
     &=&8h'(0){\rm dim}F\Omega_{4}\frac{\pi i}{5!}\Big[\frac{8+24\xi_{n}i-40\xi^{2}_{n}-24i\xi^{3}_{n}}{(\xi_{n}+i)^{4}}\Big]^{(5)}|_{\xi_{n}=i}dx'\nonumber\\
     &=&-\frac{15}{2}\pi h'(0)\Omega_{4}{\rm dim}Fdx'.
\end{eqnarray}
~\\
Similar to (3.47) and (3.48),
then we obtain
\begin{eqnarray}
&&-\frac{1}{2}f^{-1}\partial_{x_{n}}(f)
      \int_{|\xi'|=1}\int^{+\infty}_{-\infty}{\rm trace}
\Big[\pi^{+}_{\xi_{n}}\sigma_{-1}(\hat{D}_{F}^{-1})
      \times\partial^2_{\xi_{n}}\sigma_{-3}
(\hat{D}_{F}^{*}\hat{D}_{F}\hat{D}^{*}_{F})^{-1}\Big](x_0)d\xi_n\sigma(\xi')dx'\nonumber\\
      &=&(10\pi i+88\pi){\rm \Omega_{4}}dimF\cdot f^{-1}\partial_{x_{n}}(f)dx',
\end{eqnarray}
where ${\rm \Omega_{4}}$ is the canonical volume of $S_{4}.$\\
Then
\begin{eqnarray}
{\bf case~(a)~II)}&=&-\frac{15}{2}\pi h'(0)\Omega_{4}{\rm dim}Fdx'+
(10\pi i+88\pi){\rm \Omega_{4}}dimF\cdot f^{-1}\partial_{x_{n}}(f)dx',
\end{eqnarray}
~\\
 where ${\rm \Omega_{4}}$ is the canonical volume of $S_{4}.$
~\\

\noindent  {\bf case (a)~(III)}~$r=-1,l=-3,|\alpha|=j=0,k=1$.\\

By (5.2) and an integration by parts, we have

 \begin{eqnarray}
{\rm case~ (a)~(III)}&=&-\frac{1}{2}\int_{|\xi'|=1}\int^{+\infty}_{-\infty}{\rm trace} \Big[\partial_{\xi_{n}}\pi^{+}_{\xi_{n}}\sigma_{-1}(f\hat{D}_{F}^{-1})
      \times\partial_{\xi_{n}}\partial_{x_{n}}
      \sigma_{-3}
\big(f^{-1}(\hat{D}_{F}^{*})^{-1}\cdot f\hat{D}_{F}^{-1}\cdot\nonumber\\
&&f^{-1}(\hat{D}^{*}_{F})^{-1}\big)\Big](x_0)d\xi_n\sigma(\xi')dx'\nonumber\\
  &=&-\frac{1}{2}\int_{|\xi'|=1}\int^{+\infty}_{-\infty}{\rm trace}
\Big[\partial_{\xi_{n}}\pi^{+}_{\xi_{n}}\big(\sigma_{-1}(\hat{D}_{F}^{-1})\big)
      \times\partial_{\xi_{n}}\partial_{x_{n}}\sigma_{-3}
(\hat{D}_{F}^{*}\hat{D}_{F}\hat{D}^{*}_{F})^{-1}\Big](x_0)d\xi_n\sigma(\xi')dx'\nonumber\\
      &&-\frac{1}{2}f\partial_{x_{n}}(f^{-1})
      \int_{|\xi'|=1}\int^{+\infty}_{-\infty}{\rm trace}
\Big[\partial_{\xi_{n}}\pi^{+}_{\xi_{n}}\sigma_{-1}(\hat{D}_{F}^{-1})
      \times\partial_{\xi_{n}}\sigma_{-3}
(\hat{D}_{F}^{*}\hat{D}_{F}\hat{D}^{*}_{F})^{-1}\Big](x_0)d\xi_n\nonumber\\
      &&\times\sigma(\xi')dx'.
\end{eqnarray}

Similar to (3.52), (3.53) and combining (5.10), we have
\begin{eqnarray}
&&{\rm trace} \Big[\partial_{\xi_{n}}\pi^{+}_{\xi_{n}}\sigma_{-1}(\hat{D}_{F}^{-1})
      \times\partial_{\xi_{n}}\partial_{x_{n}}\sigma_{-3}
      ((\hat{D}^{*}_{F}\hat{D}_{F}\hat{D}^{*}_{F})^{-1})\Big](x_{0})|_{|\xi'|=1}\nonumber\\
&=&8h'(0){\rm dim}F\frac{8i-32\xi_{n}-8i\xi^{2}_{n}}{(\xi_{n}-i)^{5}(\xi+i)^{4}}.
\end{eqnarray}
Then
\begin{eqnarray}
&&-\frac{1}{2}\int_{|\xi'|=1}\int^{+\infty}_{-\infty}{\rm trace}
\Big[\partial_{\xi_{n}}\pi^{+}_{\xi_{n}}\big(\sigma_{-1}(\hat{D}_{F}^{-1})\big)
      \times\partial_{\xi_{n}}\partial_{x_{n}}\sigma_{-3}
(\hat{D}_{F}^{*}\hat{D}_{F}\hat{D}^{*}_{F})^{-1}\Big](x_0)d\xi_n\sigma(\xi')dx'\nonumber\\
      &=&-\frac{1}{2}\int_{|\xi'|=1}\int^{+\infty}_{-\infty}
8h'(0){\rm dim}F\frac{8i-32\xi_{n}-8i\xi^{2}_{n}}
{(\xi_{n}-i)^{5}(\xi+i)^{4}}d\xi_n\sigma(\xi')dx'\nonumber\\
     &=&-8h'(0){\rm dim}F\Omega_{4}\frac{\pi i}{4!}\Big[\frac{8i-32\xi_{n}-8i\xi^{2}_{n}}{(\xi+i)^{4}}\Big]^{(4)}|_{\xi_{n}=i}dx'\nonumber\\
     &=&\frac{25}{2}\pi h'(0)\Omega_{4}{\rm dim}Fdx',
\end{eqnarray}
and
\begin{eqnarray}
&&-\frac{1}{2}f\partial_{x_{n}}(f^{-1})
      \int_{|\xi'|=1}\int^{+\infty}_{-\infty}{\rm trace}
\Big[\partial_{\xi_{n}}\pi^{+}_{\xi_{n}}\sigma_{-1}(\hat{D}_{F}^{-1})
      \times\partial_{\xi_{n}}\sigma_{-3}
(\hat{D}_{F}^{*}\hat{D}_{F}\hat{D}^{*}_{F})^{-1}\Big](x_0) d\xi_n\sigma(\xi')dx'\nonumber\\
      &=&\frac{\pi i}{2}\cdot f\cdot\partial_{x_{n}}(f^{-1})\Omega_{4}{\rm dim}Fdx',
\end{eqnarray}
where ${\rm \Omega_{4}}$ is the canonical volume of $S_{4}.$
Then
\begin{eqnarray}
{\bf case~(a)~III)}&=&\frac{25}{2}\pi h'(0)\Omega_{4}{\rm dim}Fdx'+
\frac{\pi i}{2}\cdot f\cdot\partial_{x_{n}}(f^{-1})\Omega_{4}{\rm dim}Fdx'.
\end{eqnarray}

\noindent {\bf  case (b)}~$r=-2,l=-3,|\alpha|=j=k=0$.\\

By (5.2) and an integration by parts, we have

\begin{eqnarray}
{\rm case~ (b)}&=&-i\int_{|\xi'|=1}\int^{+\infty}_{-\infty}{\rm trace} \Big[\pi^{+}_{\xi_{n}}\sigma_{-2}(f\widehat{D}_{F}^{-1})
      \times\partial_{\xi_{n}}\sigma_{-3}\big(f^{-1}(\hat{D}_{F}^{*})^{-1}\cdot f\hat{D}_{F}^{-1}\cdot f^{-1}(\hat{D}^{*}_{F})^{-1}\big)\Big](x_0)\nonumber\\
&&\times d\xi_n\sigma(\xi')dx'\nonumber\\
&=&-i\int_{|\xi'|=1}\int^{+\infty}_{-\infty}{\rm trace} \Big[\pi^{+}_{\xi_{n}}\sigma_{-2}(\hat{D}_{F}^{-1})
      \times\partial_{\xi_{n}}
      \sigma_{-3}
\big((\hat{D}_{F}^{*}\hat{D}_{F}\hat{D}^{*}_{F})^{-1}\big)\Big]
(x_0)d\xi_n\sigma(\xi')dx'.
\end{eqnarray}

Then an application of Lemma 4.3 shows
\begin{eqnarray}
\sigma_{-2}(\hat{D}_{F}^{-1})(x_{0})&=&\frac{c(\xi)\sigma_{0}(\hat{D}_{F})(x_{0})c(\xi)}{|\xi|^{4}}+\frac{c(\xi)}{|\xi|^{6}}\sum_{j}c(dx_{j})
                    \Big[\partial_{x_{j}}(c(\xi))|\xi|^{2}-c(\xi)\partial_{x_{j}}(|\xi|^{2})\Big](x_{0})\nonumber\\
                   &=&\frac{c(\xi)\sigma_{0}(\hat{D}_{F})(x_{0})c(\xi)}{|\xi|^{4}}
                  +\frac{c(\xi)}{|\xi|^{6}}c(dx_{n})\Big[\partial_{ x_{n}}(c(\xi'))(x_{0})-c(\xi)h'(0)|\xi'|^{2}_{g^{\partial M}}\Big].
\end{eqnarray}
Hence,
\begin{equation}
\pi_{\xi_{n}}^{+}\sigma_{-2}((\hat{D}_{F})^{-1})(x_{0}):=B_{1}+B_{2}+B_{3}+B_{4},
\end{equation}
where

\begin{eqnarray}
B_1&=&\frac{-1}{4(\xi_n-i)^2}[(2+i\xi_n)c(\xi')\big(-\frac{5}{4}h'(0)c(dx_n)\big)c(\xi')
+i\xi_nc(dx_n)\big(-\frac{5}{4}h'(0)c(dx_n)\big)c(dx_n) \nonumber\\
&&+(2+i\xi_n)c(\xi')c(dx_n)\partial_{x_n}c(\xi')+ic(dx_n)\big(-\frac{5}{4}h'(0)c(dx_n)\big)c(\xi')
+ic(\xi')(-\frac{5}{4}h'(0)c(dx_n))\nonumber\\
&&\times c(dx_n)-i\partial_{x_n}c(\xi')]\nonumber\\
&=&\frac{1}{4(\xi_n-i)^2}\Big[\frac{5}{2}h'(0)c(dx_n)-\frac{5i}{2}h'(0)c(\xi')
  -(2+i\xi_n)c(\xi')c(dx_n)\partial_{\xi_n}c(\xi')+i\partial_{\xi_n}c(\xi')\Big]  ;         \\
B_2&=&-\frac{h'(0)}{2}\Big[\frac{c(dx_n)}{4i(\xi_n-i)}+\frac{c(dx_n)-ic(\xi')}{8(\xi_n-i)^2}
+\frac{3\xi_n-7i}{8(\xi_n-i)^3}[ic(\xi')-c(dx_n)]\Big].  \\
B_3&=&\frac{-1}{4(\xi_n-i)^2}\Big[(2+i\xi_n)c(\xi')pc(\xi')+i\xi_nc(dx_n)pc(dx_n)+(2+i\xi_n)c(\xi')c(dx_n)\partial_{x_n}c(\xi')\nonumber\\
&&~~~~+ic(dx_n)pc(\xi')+ic(\xi')pc(dx_n)-i\partial_{x_n}c(\xi')\Big];\\
B_4&=&\frac{-1}{4(\xi_n-i)^2}\Big[(2+i\xi_n)c(\xi')\vartheta c(\xi')+i\xi_nc(dx_n)\vartheta c(dx_n)+ic(dx_n)\vartheta c(\xi')+ic(\xi')\vartheta c(dx_n)\Big].
\end{eqnarray}

On the other hand,
\begin{equation}
\partial_{\xi_{n}}\sigma_{-3}((\hat{D}^{*}_{F}\hat{D}_{F}\hat{D}^{*}_{F})^{-1})=\frac{-4 i \xi_{n}c(\xi')}{(1+\xi_{n}^{2})^{3}}+\frac{i(1- 3\xi_{n}^{2})c(dx_{n})}
{(1+\xi_{n}^{2})^{3}}.
\end{equation}

From (5.22) and (5.26), we have
\begin{eqnarray}
&&{\rm trace }[B_1\times\partial_{\xi_n}\sigma_{-3}((\hat{D}^{*}_{F}\hat{D}_{F}\hat{D}^{*}_{F})^{-1})(x_0)]|_{|\xi'|=1}\nonumber\\
&=&{\rm tr }\Big\{ \frac{1}{4(\xi_n-i)^2}\Big[\frac{5}{2}h'(0)c(dx_n)-\frac{5i}{2}h'(0)c(\xi')
  -(2+i\xi_n)c(\xi')c(dx_n)\partial_{\xi_n}c(\xi')+i\partial_{\xi_n}c(\xi')\Big]\nonumber\\
&&\times \frac{-4i\xi_nc(\xi')+(i-3i\xi_n^{2})c(dx_n)}{(1+\xi_n^{2})^3}\Big\} \nonumber\\
&=&8h'(0)\frac{3+12i\xi_n+3\xi_n^{2}}{(\xi_n-i)^4(\xi_n+i)^3}.
\end{eqnarray}

Similarly, we obtain
\begin{eqnarray}
&&{\rm trace }[B_2\times\partial_{\xi_n}\sigma_{-3}((\hat{D}^{*}_{F}\hat{D}_{F}\hat{D}^{*}_{F})^{-1})(x_0)]|_{|\xi'|=1}\nonumber\\
&=&{\rm tr }\Big\{ -\frac{h'(0)}{2}\Big[\frac{c(dx_n)}{4i(\xi_n-i)}+\frac{c(dx_n)-ic(\xi')}{8(\xi_n-i)^2}
+\frac{3\xi_n-7i}{8(\xi_n-i)^3}[ic(\xi')-c(dx_n)]\Big]\nonumber\\
&&\times \frac{-4i\xi_nc(\xi')+(i-3i\xi_n^{2})c(dx_n)}{(1+\xi_n^{2})^3}\Big\} \nonumber\\
&=&-8h'(0)\frac{4i-11\xi_n-6i\xi_n^{2}+3\xi_n^{3}}{(\xi_n-i)^5(\xi_n+i)^3}.
\end{eqnarray}

For the signature operator case,
\begin{equation}
{\rm trace}[c(\xi')pc(\xi')c(dx_n)](x_0)={\rm trace}[pc(\xi')c(dx_n)c(\xi')](x_0)=|\xi'|^2{\rm trace}[p(x_0)c(dx_n)],
\end{equation}
 and
 \begin{eqnarray}
c(dx_n)p(x_0)&=&-\frac{1}{4}h'(0)\sum^{n-1}_{i=1}c(\tilde{e}_i)\hat{c}(\tilde{e}_i)c(\widetilde{e_n})\hat{c}(\widetilde{e_n})\nonumber\\
&=&-\frac{1}{4}h'(0)\sum^{n-1}_{i=1}[\epsilon ({\widetilde{e_i*}} )\iota
({\widetilde{e_i*}} )-\iota(\widetilde{e_i*})\epsilon(\widetilde{e_i*})][\epsilon ({\widetilde{e_n*}} )\iota
({\widetilde{e_n*}} )-\iota(\widetilde{e_n*})\epsilon(\widetilde{e_n*})].
\end{eqnarray}

 By Section 3 in \cite{Wa3}, then
 \begin{eqnarray}
&& {\rm trace}_{\wedge^m(T^*M)} \{[\epsilon ({e_i*} )\iota ({e_i*} )-\iota(e_i*)\epsilon(e_i*)]
[\epsilon ({e_n*} )\iota ({e_n*} )-\iota(e_n*)\epsilon(e_n*)]\}\nonumber\\
&=&a_{n,m}\langle e_i*,e_n* \rangle^{2}+b_{n,m}|e_i*|^2|e_n*|^2=b_{n,m},
 \end{eqnarray}

 where $b_{6,m}=\left(\begin{array}{lcr}
  \ \ 4 \\
    \  m-2
\end{array}\right)+\left(\begin{array}{lcr}
  \ \ 4 \\
    \  m
\end{array}\right)-2\left(\begin{array}{lcr}
  \ \ 4 \\
    \  m-1
\end{array}\right).$

Then
\begin{eqnarray}
{\rm tr}_{\wedge^*(T^*M)} \{[\epsilon ({\widetilde{e_i*}} )\iota ({\widetilde{e_i*}} )-\iota(\widetilde{e_i*})\epsilon(\widetilde{e_i*})]
[\epsilon ({\widetilde{e_n*}} )\iota ({\widetilde{e_n*}} )-\iota(\widetilde{e_n*})\epsilon(\widetilde{e_n*})]\}
=\sum_{m=0}^6b_{6,m}=0.
\end{eqnarray}

Hence in this case,
\begin{eqnarray}
{\rm trace}_{\wedge^*(T^*M)}[c(dx_n)p(x_0)]=0.
\end{eqnarray}

We note that $\int_{|\xi'|=1}\xi_1\cdots\xi_{2q+1}\sigma(\xi')=0$, then
${\rm trace}_{\wedge^*(T^*M)}[c(\xi')p(x_0)]$ has no contribution for computing {\bf  case (b)}.

So, we obtain
\begin{eqnarray}
&&{\rm trace }[B_3\times\partial_{\xi_n}\sigma_{-3}((\hat{D}^{*}_{F}\hat{D}_{F}\hat{D}^{*}_{F})^{-1})(x_0)]|_{|\xi'|=1}\nonumber\\
&=&{\rm trace }\Big\{\frac{-1}{4(\xi_n-i)^2}\Big[(2+i\xi_n)c(\xi')pc(\xi')+i\xi_nc(dx_n)pc(dx_n)
+(2+i\xi_n)c(\xi')c(dx_n)\partial_{x_n}c(\xi')\nonumber\\
&&~~~~+ic(dx_n)pc(\xi')+ic(\xi')pc(dx_n)-i\partial_{x_n}c(\xi')\Big] \times \frac{-4i\xi_nc(\xi')+(i-3i\xi_n^{2})c(dx_n)}{(1+\xi_n^{2})^3}\Big\} \nonumber\\
&=&8h'(0){\rm dim}F\frac{3\xi^2_n-3i\xi_n-2}{(\xi_n-i)^4(\xi_n+i)^3}.
\end{eqnarray}

Then, we have
\begin{eqnarray}
&&\text{trace}[(B_{1}+B_{2}+B_{3})\times\partial_{\xi_{n}}\sigma_{-3}((\hat{D}^{*}_{F}\hat{D}_{F}\hat{D}^{*}_{F})^{-1})](x_{0})\texttt{d}\xi_{n}\sigma(\xi')\texttt{d}x'\nonumber\\
&=&8h'(0){\rm dim}F\frac{3\xi^3_n+9i\xi^2_n+21\xi_{n}-5i}{(\xi_n-i)^5(\xi_n+i)^3}.
\end{eqnarray}

By the relation of the Clifford action and $\rm{trace}PQ=\rm{trace}QP$, then we have the equalities
\begin{eqnarray}
\rm{trace}[c(\widetilde{\emph{e}_{\emph{i}}})c(\emph{dx}_{\emph{n}})]&=&0, i<n;  \rm{trace}[c(\widetilde{\emph{e}_{\emph{i}}})c(\emph{dx}_{\emph{n}})]=-64{\rm dim}F, \emph{i}=\emph{n};\\
~~\rm{trace}[\hat{c}(\widetilde{\emph{e}_{\emph{i}}})c(\xi')]&=&\rm{trace}[\hat{c}
(\widetilde{\emph{e}_{\emph{i}}})c(\emph{dx}_{\emph{n}})]=0.
\end{eqnarray}
Then $\rm{trace}[\vartheta c(\xi')]$ has no contribution for computing {\bf  case (b)}.

Then, we have
\begin{eqnarray}
&&{\rm trace }[B_4\times\partial_{\xi_n}\sigma_{-3}((\hat{D}^{*}_{F}\hat{D}_{F}\hat{D}^{*}_{F})^{-1})]|_{|\xi'|=1}\nonumber\\
&=&{\rm trace }\Big\{\frac{-1}{4(\xi_n-i)^2}\Big[(2+i\xi_n)c(\xi')\vartheta c(\xi')+i\xi_nc(dx_n)\vartheta c(dx_n)
+ic(dx_n)\vartheta c(\xi')\nonumber\\
&&+ic(\xi')\vartheta c(dx_n)\Big]\times \frac{-4i\xi_nc(\xi')+(i-3i\xi_n^{2})c(dx_n)}{(1+\xi_n^{2})^3}\Big\} \nonumber\\
&=&\frac{i(3\xi_{n}-i)}{2(\xi_n-i)^4(\xi_n+i)^3}{\rm trace }[c(dx_{n})\vartheta]\nonumber\\
&=&-32{\rm dim}F\frac{1+3\xi_{n}i}{(\xi_n-i)^4(\xi_n+i)^3}{\rm trace }[\sigma^{F,e}_{n}].
\end{eqnarray}

From (5.35), we obtain
\begin{eqnarray}
&&-i\int_{|\xi'|=1}\int^{+\infty}_{-\infty}{\rm trace} \Big[(B_{1}+B_{2}+B_{3})
      \times\partial_{\xi_{n}}\sigma_{-3}((\hat{D}^{*}_{F}\hat{D}_{F}\hat{D}^{*}_{F})^{-1})\Big](x_0)d\xi_n\sigma(\xi')dx'\nonumber\\
&&=-8i{\rm dim}Fh'(0)\int_{|\xi'|=1}\int^{+\infty}_{-\infty}\frac{3\xi^{3}_{n}
+9\xi^{2}_{n}i+21\xi_n-5i}{(\xi_n-i)^5(\xi_n+i)^3}d\xi_n\sigma(\xi')dx' \nonumber\\
&&=-8i{\rm dim}Fh'(0)\frac{2 \pi i}{4!}\Big[\frac{3\xi^{3}_{n}+9\xi^{2}_{n}i+21\xi_n-5i}{(\xi_n+i)^3}
     \Big]^{(4)}|_{\xi_n=i}\Omega_4dx'\nonumber\\
&&=\frac{45}{2}{\rm dim}F\pi h'(0)\Omega_4dx'.
\end{eqnarray}

From (5.38), we obtain
\begin{eqnarray}
&&-i\int_{|\xi'|=1}\int^{+\infty}_{-\infty}{\rm trace} \Big[B_{4}
      \times\partial_{\xi_{n}}\sigma_{-3}((\hat{D}^{*}_{F}\hat{D}_{F}\hat{D}^{*}_{F})^{-1})\Big](x_0)d\xi_n\sigma(\xi')dx'\nonumber\\
&=&32i{\rm dim}F {\rm trace }[\sigma^{F,e}_{n}]\int_{|\xi'|=1}\int^{+\infty}_{-\infty}\frac{1+3\xi_{n}}{(\xi_n-i)^4(\xi_n+i)^3}d\xi_n\sigma(\xi')dx' \nonumber\\
&=&32i{\rm dim}F {\rm trace }[\sigma^{F,e}_{n}]\frac{2 \pi i}{3!}\Big[\frac{1+3\xi_{n}}{(\xi_n+i)^3}
     \Big]^{(4)}|_{\xi_n=i}\Omega_4dx'\nonumber\\
&=&-16{\rm dim}F {\rm trace }[\sigma^{F,e}_{n}]\Omega_4dx'.
\end{eqnarray}

Combining (5.19), (5.39) and (5.40), we have
\begin{eqnarray}
{\bf  case~b)}&=&\Big[\frac{45}{2}h'(0)-16{\rm trace }\big(\sigma^{F,e}_{n}\big)\Big]\pi {\rm dim}F \Omega_4dx'
\end{eqnarray}

 \noindent  {\bf case (c)}~$r=-1,l=-4,|\alpha|=j=k=0$.\\

By (5.2) and an integration by parts, we have
 \begin{eqnarray}
{\bf case~(c)}
&=&-i\int_{|\xi'|=1}\int^{+\infty}_{-\infty}{\rm trace} \Big[\pi^{+}_{\xi_{n}}\sigma_{-1}(f\hat{D}_{F}^{-1})
      \times\partial_{\xi_{n}}\sigma_{-4}\bigg(f^{-1}(\hat{D}_{F}^{*})^{-1}\cdot f\hat{D}_{F}^{-1}\cdot f^{-1}(\hat{D}^{*}_{F})^{-1}\bigg)\Big](x_0)\nonumber\\
      &&\times d\xi_n\sigma(\xi')dx'\nonumber\\
&=&-i\int_{|\xi'|=1}\int^{+\infty}_{-\infty}{\rm trace} \Bigg[\pi^{+}_{\xi_{n}}\sigma_{-1}(\hat{D}_{F}^{-1})
      \times\partial_{\xi_{n}}\Bigg(\sigma_{-4}(\hat{D}^{*}_{F}
      \hat{D}_{F}\hat{D}^{*}_{F})^{-1}\Bigg)\Bigg](x_0)d\xi_n\sigma(\xi')dx'\nonumber\\
&&+2if^{-1}\int_{|\xi'|=1}\int^{+\infty}_{-\infty}{\rm trace} \Bigg[\pi^{+}_{\xi_{n}}\sigma_{-1}(\hat{D}_{F}^{-1})
      \times\partial_{\xi_{n}}\Bigg(\frac{c(\xi)c(df)c(\xi)}{|\xi|^6}
      \Bigg)\Bigg](x_0)d\xi_n\sigma(\xi')dx',
\end{eqnarray}

By direct calculations, we have
\begin{equation}
\pi^{+}_{\xi_{n}}\sigma_{-1}(\hat{D}_{F}^{-1})=-\frac{c(\xi')+ic(dx_{n})}{2(\xi_{n}-i)}.
\end{equation}
In the normal coordinate, $g^{ij}(x_{0})=\delta^{j}_{i}$ and $\partial_{x_{j}}(g^{\alpha\beta})(x_{0})=0$, if $j<n$; $\partial_{x_{j}}(g^{\alpha\beta})(x_{0})=h'(0)\delta^{\alpha}_{\beta}$, if $j=n$.
So by Lemma A.2 in \cite{Wa3}, we have $\Gamma^{n}(x_{0})=\frac{5}{2}h'(0)$ and $\Gamma^{k}(x_{0})=0$ for $k<n$. By the definition of $\delta^{k}$ and Lemma 2.3 in \cite{Wa3}, we have $\delta^{n}(x_{0})=0$ and $\delta^{k}=\frac{1}{4}h'(0)c(\widetilde{e_{k}})c(\widetilde{e_{n}})$ for $k<n$. By (3.15) in \cite{WJ6}, we obtain

\begin{eqnarray}
&&\sigma_{-4}(\hat{D}^{*}_{F}\hat{D}_{F}\hat{D}^{*}_{F})^{-1}(x_{0})\nonumber\\
&=&\frac{-17-9\xi^{2}_{n}}{4(1+\xi^{2}_{n})^{4}}h'(0)c(\xi')c(dx_{n})c(\xi')+\frac{33\xi_{n}+17\xi^{3}_{n}}{2(1+\xi^{2}_{n})^{4}}h'(0)c(\xi')
+\frac{49\xi^{2}_{n}+25\xi^{4}_{n}}{2(1+\xi^{2}_{n})^{4}}h'(0)c(dx_{n})\nonumber\\
&&+\frac{1}{(1+\xi^{2}_{n})^{3}}c(\xi')c(dx_{n})\partial_{x_{n}}[c(\xi')](x_{0})-\frac{3\xi_{n}}{(1+\xi^{2}_{n})^{3}}\partial_{x_{n}}[c(\xi')](x_{0})
-\frac{2\xi_{n}}{(1+\xi^{2}_{n})^{3}}h'(0)\xi_{n}c(\xi')\nonumber\\
&&+\frac{1-\xi^{2}_{n}}{(1+\xi^{2}_{n})^{3}}h'(0)c(dx_{n})
+\frac{c(\xi)(p+\vartheta^*-\hat{c}(w^*))c(\xi)}{|\xi|^6}+\frac{\hat{c}(w)}{|\xi|^4}.
\end{eqnarray}

Then
\begin{eqnarray}
&&-i\int_{|\xi'|=1}\int^{+\infty}_{-\infty}{\rm trace} \Bigg[\pi^{+}_{\xi_{n}}\sigma_{-1}(\hat{D}_{F}^{-1})
      \times\partial_{\xi_{n}}\Bigg(\sigma_{-4}(\hat{D}^{*}_{F}
      \hat{D}_{F}\hat{D}^{*}_{F})^{-1}\Bigg)\Bigg](x_0)d\xi_n\sigma(\xi')dx'\nonumber\\
&=&-\frac{129}{2}\pi h'(0){\rm dim}F\Omega_{4}dx'+12\pi {\rm trace}\Big[\sigma_{n}^{F,e}\Big]{\rm dim}F\Omega_{4}dx'+4\pi {\rm trace}\Big[w(F,g^F)(e_{n})\Big]{\rm dim}F\Omega_{4}dx'\nonumber\\
&&-12\pi {\rm trace}\Big[w^*(F,g^F)(e_{n})\Big] {\rm dim}F\Omega_{4}dx'.
\end{eqnarray}
By $\sigma_{-1}(\hat{D}_{F}^{-1})=\sigma_{-1}(\widetilde{D}_{F}^{-1})$, similar to {\bf case~(b)} in Section 3,
and we get

\begin{eqnarray*}
&&{\rm trace} \Big[\pi^{+}_{\xi_{n}}\sigma_{-1}(\hat{D}_{F}^{-1})
      \times\partial_{\xi_{n}}\Big(\frac{c(\xi)c(df) c(\xi)}{|\xi|^{6}}\Big)\Big](x_0) \nonumber\\
&=&\frac{(4\xi_{n}i+2)i}{2(\xi_{n}+i)(1+\xi^{2}_{n})^{3}}{\rm trace}[c(\xi')c(df)]+\frac{4\xi_{n}i+2}{2(\xi_{n}+i)(1+\xi^{2}_{n})^{3}}{\rm trace}[c(dx_{n})c(df)].
\end{eqnarray*}
and
\begin{eqnarray}
&&{\rm trace} \Bigg[\pi^{+}_{\xi_{n}}\sigma_{-1}(\hat{D}_{F}^{-1})
      \times\partial_{\xi_{n}}\Big(\frac{ic(\xi)\sum\limits_j\Big[c(dx_j)|\xi|^2
+2\xi_{j}c(\xi)\Big]D_{x_j}(f^{-1})c(\xi)}{|\xi|^8}\Big)\Bigg](x_0) \nonumber\\
&=&\frac{(3\xi_{n}-i)i}{(\xi_{n}+i)(1+\xi^{2}_{n})^{4}}{\rm trace}\Bigg[c(\xi')\sum\limits_j\Big[c(dx_j)|\xi|^2
+2\xi_{j}c(\xi)\Big]D_{x_j}(f^{-1})\Bigg]\nonumber\\
&&+\frac{3\xi_{n}-i}{(\xi_{n}+i)(1+\xi^{2}_{n})^{4}}{\rm trace}\Bigg[c(dx_{n})\sum\limits_j\Big[c(dx_j)|\xi|^2
+2\xi_{j}c(\xi)\Big]D_{x_j}(f^{-1})\Bigg].
\end{eqnarray}

By the relation of the Clifford action and ${\rm trace}QP={\rm trace}PQ$, then we have the following equalities
\begin{eqnarray*}
&&{\rm trace}\Big[c(dx_{n})c(df)\Big]=-g(dx_{n},df);\\
\end{eqnarray*}
and
\begin{eqnarray*}
&&{\rm trace}\Bigg[c(dx_{n})\sum\limits_j\Big[c(dx_j)|\xi|^2
      +2\xi_{j}c(\xi)\Big]D_{x_j}(f^{-1})\Bigg]\nonumber\\
      &=&{\rm trace}\big(-{\rm id}\big)|\xi|^2\bigg(-i\partial_{x_n}(f)f^{-1}\bigg)
      +2\sum\limits_{j}\xi_j\xi_n{\rm trace}\big(-{\rm id}\big)\bigg(-i\partial_{x_j}(f)f^{-1}\bigg)\nonumber\\
      &=&-64{\rm dim}F|\xi|^2\bigg(-i\partial_{x_n}(f)f^{-1}\bigg)+2\sum\limits_{j}\xi_j\xi_n{\rm trace}\big(-{\rm id}\big)\bigg(-i\partial_{x_j}(f)f^{-1}\bigg).
\end{eqnarray*}
We note that $i<n,~\int_{|\xi'|=1}\xi_i\sigma(\xi')=0$,
so ${\rm trace}\big[c(\xi')c(df)\big]$, ${\rm trace}\Bigg[c(\xi')\sum\limits_j\Big[c(dx_j)|\xi|^2
      +2\xi_{j}c(\xi)\Big]D_{x_j}(f^{-1})\Bigg]$ and $2i\sum\limits_{j}\xi_j\xi_n\partial_{x_j}(f)f^{-1}{\rm trace}[-{\rm id}]$ have no contribution for computing {\bf case (b)}.
Then we obtain
\begin{eqnarray}
&&-2if^{-1}\int_{|\xi'|=1}\int^{+\infty}_{-\infty}{\rm trace} \Bigg[\pi^{+}_{\xi_{n}}\sigma_{-1}(\hat{D}_{F}^{-1})
      \times\partial_{\xi_{n}}\Big(\frac{c(\xi)c(df) c(\xi)}{|\xi|^{6}}\Big)\Bigg](x_0)d\xi_n\sigma(\xi')dx'\nonumber\\
      &=&\frac{3}{8f}\pi  g[dx_{n},df]\Omega_{4}dx'.
\end{eqnarray}
and
\begin{eqnarray}
&&-fi\int_{|\xi'|=1}\int^{+\infty}_{-\infty}{\rm trace} \Big[\pi^{+}_{\xi_{n}}\sigma_{-1}(\hat{D}_{F}^{-1})
      \times\partial_{\xi_{n}}\Big(\frac{ic(\xi)\sum\limits_j\Big[c(dx_j)|\xi|^2
      +2\xi_{j}c(\xi)\Big]D_{x_j}(f^{-1})c(\xi)}{|\xi|^8}\Big)\Big]\nonumber\\
      &&\times(x_0)d\xi_n\sigma(\xi')dx'\nonumber\\
      &=&-60i\partial_{x_{n}}(f)\pi  {\rm dim}F\Omega_{4}dx'.
\end{eqnarray}

Then we have
\begin{eqnarray}
{\bf case~ (c)}
&=&\Bigg\{12 {\rm trace}\Big[\sigma_{n}^{F,e}\Big]-\frac{129}{2} h'(0)
+4 {\rm trace}\Big[w(F,g^F)(e_{n})\Big]
-12 {\rm trace}\Big[w^*(F,g^F)(e_{n})\Big]\nonumber\\
&&-60i\partial_{x_{n}}(f) \Bigg\}\pi {\rm dim}F\Omega_{4}dx'+\frac{3}{8f} g(dx_{n},df)
      \pi\Omega_{4}dx'.
\end{eqnarray}

Now $\Psi$ is the sum of the {\bf case~ (a)}, {\bf case~ (b)} and {\bf case~ (c)}, then
\begin{eqnarray}
\Psi&=&\Big\{4 {\rm trace}\Big[w(F,g^F)(e_{n})\Big]-37 h'(0)-4 {\rm trace}\Big[\sigma_{n}^{F,e}\Big]-12 {\rm trace}\Big[w^*(F,g^F)(e_{n})\Big]\nonumber\\
&&+(\frac{19 i}{22f}+\frac{88}{f}-60i)\partial_{x_{n}}(f)\Big\}\pi{\rm \Omega_{4}}dimF dx'
+\frac{3}{8f} g(dx_{n},df)\pi\Omega_{4}dx'.
\end{eqnarray}

By (4.2) in \cite{Wa3}, we have
$$K=\sum_{1\leq i,j\leq n-1}K_{i,j}g^{i,j}_{\partial M};K_{i,j}=-\Gamma^{n}_{i,j},$$
and $K_{i,j}$ is the second fundamental form, or extrinsic curvature. For $n=6$, then
\begin{eqnarray}
K(x_{0})=\sum_{1\leq i,j\leq n-1}K_{i,j}(x_{0})g^{i,j}_{\partial M}(x_{0})=\sum^{5}_{i=1}K_{i,i}(x_{0})=-\frac{5}{2}h'(0).
\end{eqnarray}
Hence we conclude that
\begin{thm}
 Let M be a 6-dimensional compact manifolds   with the boundary $\partial M$. Then
 \begin{eqnarray}
&&\widetilde{Wres}\big[\pi^{+}(f\hat{D}_{F}^{-1}) \circ\pi^{+}\big(f^{-1}(\hat{D}_{F}^{*})^{-1}\cdot f\hat{D}_{F}^{-1}\cdot f^{-1}(\hat{D}^{*}_{F})^{-1}\big)\big]\nonumber\\
&=& 8\pi^{3}\int_{M}\Bigg\{\rm{trace}
 \Big[-\frac{s}{12} +\frac{3}{8}\big[\hat{c}(\omega^{*})-\hat{c}(\omega)\big]^{2}
         -\frac{1}{4}\hat{c}(\omega^{*})\hat{c}(\omega)
         -\frac{1}{4}\sum_{j}\nabla_{e_{j}}^{F}\big(\hat{c}(\omega^{*})\big)c(e_{j})\nonumber\\
       &&  +\frac{1}{4}\sum_{j}c(e_{j})\nabla_{e_{j}}^{F}\big(\hat{c}(\omega)\big)\Big]
       +4f^{-1}\Delta(f)
       +8\big\langle grad_M(f),grad_M(f^{-1})\big\rangle-5f^{-2}\Big[|grad_M(f)|^2\nonumber\\
       &&+2\Delta(f)\Big]\Bigg\}dvol_{M}+
       \int_{\partial M}\Bigg\{\Big\{4 {\rm trace}\Big[w(F,g^F)(e_{n})\Big]-4 {\rm trace}\Big[\sigma_{n}^{F,e}\Big]-12 {\rm trace}\Big[w^*(F,g^F)(e_{n})\Big]\nonumber\\
&&-37 h'(0)+\Big(\frac{19 i}{22f}+\frac{88}{f}-60i\Big)\partial_{x_{n}}(f)\Big\}{\rm dim}F
+\frac{3}{8f} g(dx_{n},df)
      \Bigg\}\pi \Omega_{4}dvol_{M}.
\end{eqnarray}
where  $s$ is the scalar curvature.
\end{thm}

\section*{Acknowledgements}
The first author is supported by DUFE202159. The
partial research of the corresponding author was supported by NSFC. 11771070.






\end{document}